\documentclass[reqno,a4paper,11pt]{article}

\usepackage{a4wide}
\usepackage[centertags]{amsmath}
\usepackage{enumerate,amsfonts,amssymb,amsthm,amsopn,cite}

\usepackage[usenames]{color}
\definecolor{citegreen}{rgb}{0,0.6,0}
\definecolor{refred}{rgb}{0.8,0,0}
\usepackage[colorlinks, citecolor=citegreen, linkcolor=refred]{hyperref}

\title{\vspace{-10mm} A Local Singularity Analysis for the Ricci Flow and its Applications to Ricci Flows with Bounded Scalar Curvature}
\author{Reto Buzano and Gianmichele Di Matteo}
\date{}

\newtheorem{theorem}{Theorem}[section]
\newtheorem{lemma}[theorem]{Lemma}
\newtheorem{corollary}[theorem]{Corollary}
\newtheorem{proposition}[theorem]{Proposition}
\newtheorem{conjecture}[theorem]{Conjecture}
\theoremstyle{definition}
\newtheorem{remark}[theorem]{Remark}

\newtheorem{definition}[theorem]{Definition}
\numberwithin{equation}{section}

\newcommand{\R}{\mathbb R}

\newcommand{\e}{\varepsilon}

\providecommand{\set}[1]{\{ {#1} \}}
\providecommand{\abs}[1]{\lvert #1\rvert}
\providecommand{\norm}[1]{\lVert #1\rVert}
\providecommand{\Norm}[1]{\Big\lVert #1\Big\rVert}

\DeclareMathOperator{\Sc}{R}
\DeclareMathOperator{\Ric}{Ric}
\DeclareMathOperator{\Rm}{Rm}

\DeclareMathOperator{\inj}{inj}

\newcommand\printaddress{{
\setlength{\parindent}{15pt}
\footnotesize~
\par
{\scshape Reto Buzano}
\newline 
Queen Mary University of London, 
School of Mathematical Sciences, 
Mile End Road, 
London E1 4NS, UK 
\newline
\emph{E-mail address:} 
\texttt{r.buzano@qmul.ac.uk}
\newline 
Universit\`a degli Studi di Torino, 
Dipartimento di Matematica,
Via Carlo Alberto 10, 
10123 Torino, Italy 
\newline
\emph{E-mail address:} 
\texttt{reto.buzano@unito.it}
\par
{\scshape Gianmichele Di Matteo}
\newline 
Queen Mary University of London, 
School of Mathematical Sciences, 
Mile End Road, 
London E1 4NS, UK
\newline
\emph{E-mail address:} 
\texttt{g.dimatteo@qmul.ac.uk}
\par
}}

\begin{document}
\maketitle
\begin{abstract}
We develop a refined singularity analysis for the Ricci flow by investigating curvature blow-up rates locally. We first introduce general definitions of Type~I and Type~II singular points and show that these are indeed the only possible types of singular points. In particular, near any singular point the Riemannian curvature tensor has to blow up at least at a Type~I rate, generalising a result of Enders, Topping and the first author that relied on a global Type~I assumption. We also prove analogous results for the Ricci tensor, as well as a localised version of Sesum's result, namely that the Ricci curvature must blow up near every singular point of a Ricci flow, again at least at a Type~I rate. Finally, we show some applications of the theory to Ricci flows with bounded scalar curvature.
\end{abstract}


\section{Introduction}
\subsection{A Refined Local Singularity Analysis}

Geometric flows typically develop singularities in finite time. A Ricci flow, that is to say a smooth one-parameter family $(M^n,g(t))$ of $n$-dimensional Riemannian manifolds satisfying $\partial_t g(t) = -2\Ric_{g(t)}$ on a time interval $t\in [0,T)$, is said to develop a singularity at time $T<\infty$ if it cannot be smoothly extended past $T$. By Hamilton's original long-time existence criterion \cite{ham1}, on closed manifolds the arrival of such a finite time singularity is characterised by the blow-up of the norm of the Riemannian curvature tensor. This result was improved by Sesum \cite{ses} who showed that if the Ricci curvature remains bounded on $[0,T)$ then the flow can be extended past $T$. Hamilton's and Sesum's results have then been generalised to extension theorems under a wide variety of other pointwise or integral curvature bounds on closed manifolds or complete manifolds with bounded curvature, see for example \cite{ct15,che,dim,he14,kno,kot,wan1,ye,ye1} for a non-exhaustive list.

Many of the results currently present in the literature study and classify singularities from a \emph{global} point of view, without considering where the flow becomes singular. Of particular importance is the work of Hamilton, who introduced different notions of finite time singularities, called Type~I and Type~II singularities, distinguishing them by the rate at which the \emph{maximal curvature} blows up at the singular time, see \cite{ham2}. Motivated by his theory, the results developed by Angenent and Knopf in \cite{ang}, as well as several recent examples of Type~II singularities, e.g. \cite{ang2,app,dig,sto,wu}, we introduce a \emph{local} analysis of singularities and curvature blow-up rates.

As mentioned above, it is well known that a closed Ricci flow $(M,g(t))_{t\in[0,T)}$ cannot be smoothly extended past time $T<\infty$ if and only if the Riemannian curvature tensor $\Rm$ blows up, i.e. 
\begin{equation}\label{eq.curvatureblowup}
\limsup_{t\nearrow T}\,\sup_M\abs{\Rm(\cdot,t)}_{g(t)}=\infty
\end{equation}
or equivalently
\begin{equation}\label{eq.typeIlower}
(T-t) \sup_M\abs{\Rm(\cdot,t)}_{g(t)}\geq 1/8, \quad \forall t\in[0,T).
\end{equation}
Property \eqref{eq.typeIlower} follows from \eqref{eq.curvatureblowup} by a maximum principle argument applied to the evolution equation of the Riemannian curvature tensor along the Ricci flow. In the noncompact case, there exist Ricci flows satisfying (\ref{eq.curvatureblowup}) that can be smoothly extended past $T$ as well as flows that become singular at time $T$ but have unbounded curvatures for some or all earlier times, see \cite{cab,gie} for instance. In particular, if the Ricci flow is not closed, (\ref{eq.curvatureblowup}) is necessary but not sufficient for the flow to develop a singularity. In this work however, we (almost) always restrict ourselves to flows that have bounded curvature on every compact subinterval of $[0,T)$, i.e.
\begin{equation}\label{eq.bddcurv}
\forall t\in[0,T)\, \exists C<\infty \text{ such that } \abs{\Rm} \leq C \text{ on } M \times [0,t],
\end{equation}
so that $T$ is the first time of curvature blow-up. For complete Ricci flows satisfying \eqref{eq.bddcurv}, the conditions \eqref{eq.curvatureblowup} and \eqref{eq.typeIlower} are still equivalent. The reader may be misled to believe that assuming \eqref{eq.curvatureblowup}--\eqref{eq.bddcurv} forces the existence of singular points, but we notice that in \cite{car}, the authors constructed Ricci flows satisfying these assumptions and becoming singular only at spatial infinity. 

A Ricci flow on $[0,T)$ satisfying \eqref{eq.curvatureblowup} is said to be of \emph{Type~I} if there exists an upper bound analogous to \eqref{eq.typeIlower}, in other words, if there exists a constant $C$ such that
\begin{equation}\label{eq.typeIupper}
(T-t)\sup_M\abs{\Rm(\cdot,t)}_{g(t)}\leq C, \quad \forall t\in[0,T).
\end{equation}
If no such $C$ exists, meaning that
\begin{equation}\label{eq.typeIIglobal}
\limsup_{t\nearrow T}\,(T-t)\sup_M\abs{\Rm(\cdot,t)}_{g(t)} = \infty,
\end{equation}
the Ricci flow is said to be of \emph{Type~II}. In our first definition we localise these concepts.

\begin{definition}[Type~I and Type~II Singular Points]\label{def.singpoints}
Let $(M,g(t))$ be a Ricci flow on $[0,T)$, $T<\infty$, satisfying \eqref{eq.curvatureblowup} and \eqref{eq.bddcurv}. For any fixed $t\in[0,T)$, we consider the parabolically rescaled Ricci flow $\widetilde{g}_t(s):= (T-t)^{-1}g(t+(T-t)s)$ defined for $s\in[-\frac{t}{T-t},1)$.
\begin{enumerate}[i)]
\item We say that a point $p \in M$ is a \emph{singular point} if for any neighbourhood $U$ of $p$, the Riemannian curvature becomes unbounded on $U$ as $t$ approaches $T$. The \emph{singular set} $\Sigma$ is the set of all such points and the \emph{regular set} $\mathfrak{Reg}$ consists of the complement of $\Sigma$.
\item We say that a point $p \in M$ is a \emph{Type~I singular point} if there exist constants $c_I,C_I, r_I>0$ such that we have
\begin{equation}\label{eq.typeI}
c_I < \limsup_{t\nearrow T} \sup_{B_{\widetilde{g}_t(0)}(p,r_I)\times(-r_I^2,r_I^2)} \,\abs{\Rm_{\widetilde{g}_t}}_{\widetilde{g}_t} \leq C_I.
\end{equation}
We denote the set of such points by $\Sigma_I$ and call it the \emph{Type~I singular set}.
\item We say that a point $p$ is a \emph{Type~II singular point} if for any $r>0$ we have
\begin{equation}\label{eq.typeII}
\limsup_{t\nearrow T} \sup_{B_{\widetilde{g}_t(0)}(p,r)\times(-r^2,r^2)} \, \abs{\Rm_{\widetilde{g}_t}}_{\widetilde{g}_t}=\infty.
\end{equation}
We denote the set of such points by $\Sigma_{II}$ and call it the \emph{Type~II singular set}.
\end{enumerate}
\end{definition}

A somewhat related local definition of Type~I and Type~II singular points for the mean curvature flow has appeared in a very recent preprint \cite{SS20}, but it differs from ours in the sense that it uses backwards parabolic cylinders based at the singular time, while we use forwards and backwards parabolic cylinders based at regular times. This subtle difference turns out to be crucial for our results below.

Part i) of Definition \ref{def.singpoints} implies that if $p\in\Sigma$ is a singular point, then there exists a sequence of space-time points $(p_i,t_i)$ such that $p_i \to p$, $t_i \to T$ and $\abs{\Rm}(p_i,t_i)\to \infty$. We call such a sequence an \emph{essential blow-up sequence}. Parts ii) and iii) do not only consider the rate of curvature blow-up, but also take into account the rate of convergence of essential blow-up sequences to the singular point. In Section \ref{sec.pointwise}, we give a heuristic explanation for the precise choices of $\Sigma_I$ and $\Sigma_{II}$.

Note that if the Ricci flow in consideration is globally of Type~I in the sense of \eqref{eq.typeIupper}, then so are all the rescaled flows $(M,\widetilde{g}_t)$ -- with the same constant $C$ -- and therefore we immediately obtain the upper bound in \eqref{eq.typeI} for any radius $r_I<1$ with $C_I = \frac{C}{1-r_I^2}$. The first author, together with Enders and Topping, showed that under a global Type~I assumption the Riemannian curvature blows up at a Type~I rate at any singular point, i.e. we also obtain the lower bound in \eqref{eq.typeI}, see Theorem $3.2$ in \cite{end1}. We extend this result to great generality in the following theorem.

\begin{theorem}[Decomposition of the Singular Set]\label{thm.singdecomposition}
Let $(M,g(t))$ be a Ricci flow on $[0,T)$, $T<\infty$, satisfying \eqref{eq.curvatureblowup} and \eqref{eq.bddcurv}. Then $\Sigma=\Sigma_I \cup \Sigma_{II}$.
\end{theorem}
We point out that this theorem not only shows that for every singular point $p\in\Sigma$ there is an essential blow-up sequence $(p_i,t_i)$ such that the curvature blows up with
\begin{equation}\label{eq.curvblowupalongseq}
\lim_{i\to\infty}\, (T-t_i)\abs{\Rm(p_i,t_i)}_{g(t_i)} >0,
\end{equation}
but this blow-up sequence can be chosen to satisfy also
\begin{equation}\label{eq.convergenceofseq}
\limsup_{i\to\infty}\, \frac{d^2_{g(t_i)}(p_i,p)}{T-t_i} <\infty.
\end{equation}
In fact, if $p$ is a Type~II singular point, then there exists an essential blow-up sequence $(p_i,t_i)$ with
\begin{equation}\label{eq.essseqTypeII}
\lim_{i\to\infty}\, (T-t_i)\abs{\Rm(p_i,t_i)}_{g(t_i)} =\infty \qquad\text{and}\qquad \lim_{i\to\infty}\, \frac{d^2_{g(t_i)}(p_i,p)}{T-t_i} =0.
\end{equation}
While Theorem $3.2$ in \cite{end1} relies on the fact that Type~I blow-up limits are nontrivial gradient shrinking solitons and on Perelman's pseudolocality theorem from \cite{per}, our Theorem \ref{thm.singdecomposition} is proved using much more elementary estimates. The main technical tool is the following notion of Riemann (or regularity) scale. 

\begin{definition}[Parabolic Cylinders and Riemann Scales]\label{def.riemannscale}
Let$(M,g(t))$ be a Ricci flow defined on $[0,T)$ and let $(p,t) \in M \times [0,T)$ be a space-time point.
\begin{enumerate}[i)]
\item For $r>0$, we define the \emph{parabolic cylinder} $\mathcal{P}(p,t,r)$ with centre $(p,t)$ and radius $r$ by
\begin{equation}\label{eq.paraboliccylinderdef}
\mathcal{P}(p,t,r) := B_{g(t)}(p,r)\times(\max \{ t-r^2, 0\},\min \{t+r^2,T\}).
\end{equation}
\item We define the \emph{Riemann scale} $r_{\Rm}(p,t)$ at (p,t) by
\begin{equation}\label{eq.Rmscalefb}
r_{\Rm}(p,t) := \sup \{ r>0 \mid \abs{\Rm}<r^{-2} \text{ on } \mathcal{P}(p,t,r) \}.
\end{equation}
If $(M,g(t))$ is flat for every $t\in[0,T)$, we set $r_{\Rm}(p,t)=+\infty$. Moreover, by slight abuse of notation, we may sometimes write $\mathcal{P}(p,t,r_{\Rm})$ for $\mathcal{P}(p,t,r_{\Rm}(p,t))$.
\item We define the \emph{time-slice Riemann scale} $\widetilde{r}_{\Rm}(p,t)$ at $(p,t)$ by
\begin{equation}\label{eq.Rmscaletimeslice}
\widetilde{r}_{\Rm}(p,t):= \sup \set{r>0 \mid \abs{\Rm}<r^{-2} \text{ on } B_{g(t)}(p,r)}.
\end{equation}
When the flow is flat at time $t$, we set $\widetilde{r}_{\Rm}(p,t)=+\infty$. Clearly $\widetilde{r}_{\Rm}(p,t) \geq r_{\Rm}(p,t)$.
\end{enumerate}
\end{definition}

Notions of Riemann scales similar to the ones above have first been defined for static manifolds (see e.g.~\cite{and}) and both \eqref{eq.Rmscaletimeslice} as well as a definition involving backwards parabolic cylinders have appeared in various results about the Ricci flow (see e.g.~\cite{bam,bamz,hei}), but our definition \eqref{eq.Rmscalefb} using \emph{forwards and backwards} parabolic cylinders seems to have some advantages. In particular, we can prove that this Riemann scale is Lipschitz continuous in space and H\"older continuous in time, see Theorem \ref{thm.lipschitzholder}, a result which has several interesting corollaries, for example a local Harnack-type inequality (Corollary \ref{harnackcor}), estimating the infimum and supremum of $r_{\Rm}$ on a smaller parabolic cylinder $\mathcal{P}(p,t,a_1\,r_{\Rm}(p,t))$, with $a_1\in(0,1)$, by its value at the center of the cylinder. A similar definition of regularity scale using forwards and backwards cylinders has previously been given for mean curvature flow in \cite{chg2}.

Using the Riemann scale, we can give an alternative characterisation of the different types of singular points which should be compared to their global counterparts in \eqref{eq.curvatureblowup}--\eqref{eq.typeIIglobal}.

\begin{theorem}[Alternative Characterisation of Singular Sets]\label{thm.singRiemannscale}
Let $(M,g(t))$ be a Ricci flow on $[0,T)$, $T<\infty$, satisfying \eqref{eq.curvatureblowup} and \eqref{eq.bddcurv}. Let $\Sigma$, $\Sigma_I$, and $\Sigma_{II}$ be given by Definition \ref{def.singpoints}. Then
\begin{enumerate}[i)]
\item $p\in \Sigma$ if and only if $\limsup_{t\nearrow T}\, r_{\Rm}^{-2}(p,t)=\infty$.
\item $p\in\Sigma_I$ if and only if for some $0<\widetilde{c}_I, \widetilde{C}_I$ we have $\widetilde{c}_I < \limsup_{t\nearrow T}\, (T-t) r_{\Rm}^{-2}(p,t) \leq \widetilde{C}_I$.
\item $p\in\Sigma_{II}$ if and only if $\limsup_{t\nearrow T}\, (T-t) r_{\Rm}^{-2}(p,t)=\infty$.
\end{enumerate}
\end{theorem}

We also note that the Riemann scale cannot oscillate between the Type~I rate and a lower rate in the sense that if $p\in\Sigma_I$ then we also obtain $1 \leq \liminf_{t\nearrow T}\, (T-t) r_{\Rm}^{-2}(p,t)$. This is basically a feature of the definition in \eqref{eq.Rmscalefb}. We expect that a similar result should also be true for Type~II singular points and thus conjecture that in Theorem \ref{thm.singRiemannscale}, we can replace each instance of $\limsup_{t\nearrow T}$ with $\liminf_{t\nearrow T}$. We provide some evidence for this at the end of Section \ref{sec.pointwise}.

\subsection{An Integral Concentration Result and a Density Function}
Next, we study an integral characterisation of the different types of singular points. In \cite{dim}, the second author studied space-time integral curvature bounds along the Ricci flow. A key role in his analysis is given by the following concept.
\begin{definition}[Optimal Pair]\label{def.optimal}
A pair $(\alpha,\beta) \in (1,\infty)\times (1,\infty)$ is said to be \emph{optimal} if
\begin{equation}
 \alpha= \frac{n}{2} \frac{\beta}{\beta-1}.
\end{equation}
\end{definition}
In particular, he showed that a Ricci flow can always be extended if the space-time integral norm 
\begin{equation*}
\norm{\Rm}_{\alpha,\beta,M\times[0,T)} :=\bigg(\int_0^T \bigg(\int_M \abs{\Rm}^\alpha d\mu_s\bigg)^{\beta/\alpha} ds\bigg)^{1/\beta}
\end{equation*}
is finite, for some optimal pair $(\alpha,\beta)$. Heuristically, one would then expect that this integral norm of the curvature should concentrate in neighbourhoods of a singular point. Such a result would in theory be equivalent to a Harnack inequality for $\abs{\Rm}$ on parabolic cylinders near the singular point and would in particular guarantee that $\abs{\Rm}$ blows up \emph{at} (rather than \emph{near}) every singular point, a result that seems too hard to achieve with the tools developed here. Moreover, for regular points one must clearly weaken any such claim: in fact, considering flat points in the Ricci flow starting at an immersed two-torus, we see that they are surrounded by non-flat points, so no Harnack inequality in any parabolic cylinder centred at them can hold without a correction term. Hence, we instead consider space-time integral norms of $r_{\Rm}^{-2}$ on parabolic cylinders $\mathcal{P}(p,t,a_1\,r_{\Rm})$ on which our Harnack-type inequality for the Riemann scale holds, obtaining the following $\e$-regularity result: if $p\in\Sigma$, then for $t$ sufficiently close to $T$, we have 
\begin{equation}\label{eq.epsreg1}
0 < C_2 \leq \norm{r_{\Rm}^{-2}}_{\alpha,\beta,\mathcal{P}(p,t,a_1\,r_{\Rm})} \leq C_3 <\infty,
\end{equation}
where $a_1\in(0,1)$ is the constant from Corollary \ref{harnackcor}. See Theorem \ref{riemannscaleintegralbounds} for the precise statement and the dependence of the constants $C_2$ and $C_3$. In order to distinguish between the different types of singular points, we then compare their Riemann scale with the Type~I rate, so we consider the following.
\begin{definition}[Singular Density]\label{def.density}
Given a Ricci flow $(M,g(t))$ on $[0,T)$, $T<\infty$, satisfying \eqref{eq.curvatureblowup} and \eqref{eq.bddcurv}, as well as an optimal pair $(\alpha,\beta)$, we define \emph{the singular density function at time $T$} of the flow to be the function $\Theta \colon M \to [0,+\infty]$ given by
\begin{equation}
\Theta(p) := \liminf_{t \nearrow T} \Norm{\frac{1}{T-s}}_{\alpha,\beta,\mathcal{P}(p,t,a_1\,r_{\Rm})}
\end{equation}
where $a_1\in(0,1)$ is again the constant from the Harnack-type result in Corollary \ref{harnackcor}.
\end{definition}
It is then easy to prove the following alternative characterisation of the singular sets.
\begin{theorem}[Integral Classification of the Singular Sets] \label{thm.singintegralcharact}
Let $(M,g(t))$ be a complete Ricci flow on $[0,T)$, $T<\infty$, satisfying \eqref{eq.curvatureblowup} and \eqref{eq.bddcurv} as well as $\inj(M,g(0))>0$. Let $(\alpha,\beta)$ be an optimal pair and $\Theta$ the associated singular density function as defined above. Then we have $\Sigma_I =\set{p\in M \mid \Theta(p) \in (0,\infty)}$ and $\Sigma_{II} =\set{p\in M \mid \Theta(p) =0}$.
\end{theorem}

We want to remark that one of the main reasons for introducing an integral concept of a density function is that we expect the following conjecture to be true.
\begin{conjecture}
The density $\Theta$ is lower semi continuous with respect to the topology of $(M,g(t))$, and therefore, in particular, its zero-level set $\Sigma_{II}$ is closed.
\end{conjecture}

\subsection{The Ricci Singular Sets and a Localised Version of Sesum's Result}

In Definition \ref{def.ricciscale}, we give definitions of Ricci scale $r_{\Ric}(p,t)$ and time-slice Ricci scale $\widetilde{r}_{\Ric}(p,t)$ similar to Definition \ref{def.riemannscale} above. In a similar way, in Section \ref{sec.Ricci}, we then also introduce versions of Definition \ref{def.singpoints} involving the Ricci curvature tensor instead of the full Riemannian curvature tensor, i.e. we define the \emph{Ricci singular set} $\Sigma^{\Ric}$ as well as the \emph{Type~I} and \emph{Type~II Ricci singular sets} $\Sigma_{I}^{\Ric}$ and $\Sigma_{II}^{\Ric}$ in Definition \ref{def.riccisingpoints}. As it turns out, the theorems above all have a direct Ricci curvature counterpart: in particular, we prove the Lipschitz-H\"older continuity of the Ricci scale in Theorem \ref{thm.riccilipschitzholder}, the alternative characterisation of singular sets in Theorem \ref{thm.singRicciscale}, and finally the decomposition of the Ricci singular set $\Sigma^{\Ric} = \Sigma_{I}^{\Ric} \cup \Sigma_{II}^{\Ric}$ in Corollary \ref{cor.riccisingdecomposition}. We also obtain an $\e$-regularity similar to \eqref{eq.epsreg1}, see Theorem \ref{ricciscaleintegralbounds}. 

All of these results follow rather similarly to their Riemann scale versions. The main new result of Section \ref{sec.Ricci} instead is a localised version of the result that the Ricci curvature blows up at a singularity of the Ricci flow. In \cite{ses}, Sesum proved that if $(M,g(t))$ is a closed Ricci flow maximally defined on $[0,T)$, the Ricci curvature tensor $\Ric$ satisfies 
\begin{equation}\label{eq.riccicurvatureblowup}
\limsup_{t\nearrow T}\,\sup_M\abs{\Ric(\cdot,t)}_{g(t)}=\infty.
\end{equation}
This was later extended by Ma-Cheng to complete and bounded curvature Ricci flows. In \cite{wan}, Wang strengthened this result, showing that, similar to \eqref{eq.typeIlower}, if $T$ is the singular time of a closed Ricci flow $(M,g(t))$, one has
\begin{equation}\label{eq.riccitypeIlower}
(T-t) \sup_M\abs{\Ric(\cdot,t)}_{g(t)}\geq \eta_1, \quad \forall t\in[0,T).
\end{equation}
Here $\eta_1=\eta_1(n,\kappa)$ is a constant depending on the dimension of $M$ and the non-collapsing constant $\kappa$ of the flow. We generalise these results to the local setting, showing that any singular point is also a Ricci singular point -- the other direction is obviously true -- hence obtaining the following theorem.
\begin{theorem}[Singular Points are Ricci Singular Points]\label{thm.bridge}
Let $(M,g(t))$ be a complete Ricci flow on $[0,T)$, $T<\infty$, satisfying \eqref{eq.curvatureblowup} and \eqref{eq.bddcurv}. Suppose that the initial slice satisfies $\inj(M,g(0))>0$. Then $\Sigma=\Sigma^{\Ric}$.
\end{theorem}

We remark that this theorem is not a direct consequence of the local curvature bounds obtained in \cite{che}, \cite{wan}, and \cite{kot}. While their results bound \emph{the oscillations} of $\Rm$ (locally), we instead require a bound on the absolute value. Our proof relies both on the original ideas of Sesum \cite{ses} as well as on the characterisation of Ricci singular points in Theorem \ref{thm.singRicciscale}. A direct corollary of this result is the following.

\begin{corollary}\label{cor.bridge}
Let $(M,g(t))$ be a complete Ricci flow on $[0,T)$, $T<\infty$, satisfying \eqref{eq.curvatureblowup} and \eqref{eq.bddcurv}. Suppose that the initial slice $(M,g(0))$ satisfies $\inj(M,g(0))>0$. Then we have the inclusions $\Sigma_I \subseteq \Sigma_I^{\Ric}$ and $\Sigma_{II} \supseteq \Sigma_{II}^{\Ric}$ as well as the identity 
\begin{equation*}
\Sigma_{II} \setminus \Sigma_{II}^{\Ric}=\Sigma_I^{\Ric} \setminus \Sigma_I.
\end{equation*}
\end{corollary}

\subsection{Applications to Bounded Scalar Curvature Ricci Flows}

In the final section we discuss how and to which extent the local theory described above can be applied to the study of bounded scalar curvature Ricci flows. It has been conjectured that a bound on the scalar curvature could potentially also be sufficient to extend the flow. In dimension three, this is a consequence of the Hamilton-Ivey pinching estimate \cite{ham2,I93} while in higher dimensions it is known to be true for Type~I Ricci flows by Enders, Topping and the first author \cite{end1} as well as in the K\"ahler case by Zhang \cite{zhang}. In recent years, this conjecture has been the focus of many interesting new developments, see for example \cite{bam,bamz,chwan,chwan1,chwan2,les1,sim,sim1,wan,z12} and the references therein, but without the Type~I or K\"ahler assumption the conjecture still remains open in dimensions $n\geq 4$. 

In order for our theory to apply, we need to exclude badly behaved singular points, at which the Ricci curvature blows up at a lower rate than their Ricci curvature scale. We will show that well behaved singularities cannot occur in dimensions lower than eight and that in higher dimensions the well-behaved singular set has codimension at least eight. 

In the context of Ricci flow, the lack of an ambient space with respect to which we can measure this dimension forces us to consider dimensional bounds in terms of volume estimates on the singular set. In fact, the dimension of the singular set is related to the rate of convergence of its volume to zero as $t$ approaches the singular time. This approach revealed useful in the study of Type I Ricci flows by Gianniotis (see \cite{gia1,gia2}) and we briefly recall the heuristic behind it. An actual estimate on the (intrinsic) Minkowski content cannot be available in general. Indeed, if we consider the Ricci flow of a round sphere $\mathbb{S}^n$, we see that the singular set coincides with the entire manifold, so for every time $t$ the singular set is $n$-dimensional. On the other hand, the flow collapses the sphere to a single point as $t$ approaches the final time, and one can easily see that the volume $\mu_{g(t)}(\mathbb{S}^n) \sim (\sqrt{T-t})^n$ as $t \rightarrow T$, so that the volume of the singular set goes to zero at the fastest possible rate, which means that it may be interpreted as (Minkowski) $0$-dimensional. We phrase our codimension eight result in the sense of such a decay estimate.

It is worth mentioning that addressing the issue from an extrinsic point of view is in principle also possible: we could study bounds on the dimension either in the space-time structure developed by Kleiner and Lott in \cite{KL}, or in the final time slice of the flow, which can be endowed with a pseudo-metric structure by the works of Bamler-Zhang \cite{bamz2} (one might want to pass to the quotient metric space). Finally, we could also follow the work of Bamler \cite{bam} to pass to a singular limit and estimate the singular set there. Our decay estimates should essentially be equivalent to this last approach.

As mentioned above, for these results we need to exclude badly behaved singular points. This technical assumption allows us to compare the Ricci scale to the square root of the Riemann scale at the points in consideration. It would be interesting to remove this assumption or in fact to rule out such badly behaved points not only in the bounded scalar curvature case but possibly even for general Ricci flows. To phrase our results precisely, let us consider a Ricci flow $(M,g(t))$ defined on $[0,T)$, and let us assume that $\abs{\Ric}$ is not identically $0$. For every $\delta \in(0,1)$, we can consider the set of $\delta$-well behaved points
\begin{equation}\label{eq.wellbehaved}
G_\delta = G_{\delta,t_1}:= \big\{ q \in M \mid \delta a_0 r^{-2}_{\Ric}(q,t) < \abs{\Ric}(q,t), \text{ for all }t\in[t_1,T)\big\}.
\end{equation}
(The constant $a_0=\sqrt{n}(n-1)$ is added for convenience, since it implies $r_{\Ric}\geq r_{\Rm}$, compare with Remark \ref{rem.notourdefinition}, and the value of $t_1$ will be chosen suitably below.) We point out that points in these sets are well behaved in the above mentioned sense \emph{in a uniform way} near the singular time $T$ and obviously the size of these sets increases as $\delta$ goes to $0$. Set $\Sigma_{\delta} := \Sigma \cap G_\delta$. Implicit in our method of proof is the fact that $\Sigma_{\delta} \subseteq \Sigma_{II}$ for any $\delta>0$. It is not clear whether we also have
\begin{equation}
\Sigma_{II}= \overline{\bigcup_{\delta \in (0,1)} \Sigma_{\delta}}.
\end{equation}
Similarly, we define well-behaved blow-up sequences as follows: A sequence $(p_i,t_i)$ with $t_i\nearrow T$ and $r_{\Ric}(p_i,t_i)\to 0$ is said to be $\delta$-well behaved, if for sufficiently large $i$, the $\sqrt{\delta} \widetilde{r}_{\Ric}(p_i,t_i)$-ball around $(p_i,t_i)$ contains only $\delta$-well behaved points, that is for all $q \in B_{g(t_i)}(p_i,t_i,\sqrt{\delta} \widetilde{r}_{\Ric}(p_i,t_i))$ we have $\delta a_0 r^{-2}_{\Ric}(q,t_i) < \abs{\Ric}(q,t_i)$.

Our first result is a non-existence result of well-behaved singularities in dimensions $n<8$.

\begin{theorem}[No Well-Behaved Singularities in Dimensions $n<8$]\label{th.nosingularities}
Let $(M,g(t))$ be a Ricci flow on a closed manifold $M$ of dimension $n<8$, defined on $[0,T)$, $T<+\infty$. Assume that the scalar curvature is uniformly bounded, $\abs{\Sc} \le n(n-1)R_0<\infty$ on $M \times [0,T)$. Then for any $t_1 \in (0,T)$ and $\delta \in (0,1)$ there cannot be any $\delta$-well-behaved blow-up sequences. Moreover, if $M=G_\delta$ for some $\delta>0$, then the flow can be smoothly extended past time $T$.
\end{theorem}

As a corollary, we obtain an extension result under the slightly stronger assumption of an injectivity radius bound.

\begin{corollary}[No Singularities in Dimensions $n<8$ Under an Injectivity Radius Bound]\label{cor.nosingularities}
Let $(M,g(t))$ be a Ricci flow on a closed manifold $M$ of dimension $n<8$, defined on $[0,T)$, $T<+\infty$. Assume that the scalar curvature is uniformly bounded, $\abs{\Sc} \le n(n-1)R_0<\infty$ on $M \times [0,T)$ and the injectivity radius is bounded from below by 
\begin{equation}\label{eq.injbound}
\inj(M,g(t)) \geq \alpha \Big(\sup_{M\times[0,t]} \abs{\Ric}\Big)^{-1/2}
\end{equation}
for some $\alpha>0$. Then the flow can be smoothly extended past time $T$.
\end{corollary}

In our last result, we give a codimension estimate for the well-behaved singular set.

\begin{theorem}[The Well-Behaved Singular Set has Codimension $8$]\label{th.codimension8}
For any $n \in \mathbb{N}$, $R_0>0$, $i_0>0$, $k_0>0$, $\delta\in(0,1)$ and $d \in (0,8)$ there exist a constant $E=E(n,R_0,i_0,k_0,T,\delta,d)$ and a time $t_1=t_1(n,i_0,k_0,R_0,T,\delta)\in(0,T)$ such that the following statement holds. Let $(M,g(t))$ be a Ricci flow on a closed manifold $M$ of dimension $n$, defined on $[0,T)$, $T<+\infty$. Assume that the scalar curvature is uniformly bounded, $\abs{\Sc} \le n(n-1)R_0<\infty$ on $M \times [0,T)$,  and that the initial metric satisfies $\inj(M,g(0))>i_0$ and $\Ric_{g(0)}\ge -(n-1)k_0 g(0)$.

Set $\Sigma_\delta = \Sigma \cap G_{\delta,t_1}$. Then for any $p \in M$ and $t \in [t_1,T)$ we have 
\begin{equation}\label{eq.codimest}
\mu_{g(t)} \big(\Sigma_{\delta} \cap B_{g(t)}(p,\tfrac{1}{2})\big) \le E (\sqrt{T-t})^d.
\end{equation}
\end{theorem}
The constant $E$ in Theorem \ref{th.codimension8} degenerates as $d$ approaches $8$ or $\delta$ approaches $0$. This is due to our application of Proposition $6.4$ of \cite{bam}, which plays a crucial role in our argument. Preventing this degeneration for $d=8$ would correspond in this context to the finiteness of the $(n-8)$-dimensional measure of $\Sigma_{\delta}$.

The importance of considering well-behaved points in $G_\delta$ consists on the fact that they verify a property (see Theorem \ref{thm.squareroot}) which is one of the key ingredients of our proof of Theorem \ref{th.codimension8}, namely we can compare the \emph{square} of the parabolic Ricci scale with the parabolic Riemann scale ($r_{\Ric}^{2}\gtrsim r_{\Rm}$) at these points. We note that such an estimate also implies the corresponding estimate $\widetilde{r}_{\Ric}^{2}\gtrsim \widetilde{r}_{\Rm}$ for the time-slice scales.

Let us briefly describe the proofs of these results. The proof of Theorem \ref{th.nosingularities} is via an argument by contradiction. If the flow contains a $\delta$-well-behaved blow-up sequence, then by the powerful integral bound of Theorem $1.7$ in \cite{bam} by Bamler (see also \cite{sim} for a similar result in dimension four), and the estimate $\widetilde{r}_{\Ric}^{\,2}\gtrsim \widetilde{r}_{\Rm}$, we see that $\widetilde{r}_{\Ric}^{\,-1}$ has infinitesimal $L^{8-4\e}$-norm along the sequence $(p_i,t_i)$ for $\e>0$ small enough. On the other hand this norm is bounded away from zero by our Ricci Scale Concentration Lemma \ref{lemma.L8lowerbound}, yielding the desired contradiction for the first statement. To obtain the second statement, assume that $M=G_\delta$ for some $\delta>0$. If the flow develops a singularity, we can pick a sequence of space-time points $(p_i,t_i)$ along which the Ricci curvature blows up by Sesum's result \cite{ses}, in particular $r_{\Ric}(p_i,t_i)\to 0$. Since $M=G_\delta$, this blow-up sequence must be $\delta$-well behaved. But such a sequence cannot exist by the first part of the theorem.

The main idea to prove the codimension eight estimate is to use our localised version of Sesum's result from Theorem \ref{thm.bridge} and to apply suitably Proposition $6.4$ of \cite{bam}, which gives a bound on the volumes of lower level sets of the time-slice Riemann scale (or equivalently of the parabolic Riemann scale). Theorem \ref{thm.bridge} ensures that the singular set is contained in small lower level sets of the Ricci scale. Since $r_{\Ric}^{2}\gtrsim r_{\Rm}$, these level sets are comparable to drastically smaller lower level sets of the Riemann scale yielding a significant improvement in the volume bound (for well behaved points). It is worth remarking that a straightforward application of Proposition $6.4$ of \cite{bam} would give a codimension $4$ result on the entire singular set $\Sigma$ in the sense of \eqref{eq.codimest}.

We conclude the introduction outlining the structure of the paper. In Section \ref{sec.pointwise} we carry out the pointwise analysis of the Riemann scale described above. The results involving mixed integral norms will then be proven in Section \ref{sec.mixedintegral}. The results involving the Ricci scale can be found in Section \ref{sec.Ricci}. Finally, in Section \ref{sec.scalar2}, we prove our result about bounded scalar curvature flows, Theorems \ref{th.nosingularities} and \ref{th.codimension8}.

\subsection{Acknowledgements.} 
We would like to thank Andrea Malchiodi, Miles Simon and Peter Topping for their interest, comments and several stimulating discussions. We also thank an anonymous referee for some insightful remarks that improved our presentation. The first author has been supported by the EPSRC grant EP/S012907/1.


\section{Pointwise Analysis of the Singular Sets}\label{sec.pointwise}
Let us start this section with some heuristic ideas behind Definition \ref{def.singpoints} and a comparison to existing results in the literature. In \cite{end1}, the first author together with Enders and Topping gave an alternative definition of the Type~I singular set $\Sigma_{I}$ as the set of points $p$ for which there exists an essential blow-up sequence $(p_i,t_i)$ satisfying \eqref{eq.curvblowupalongseq}, see Definition $1.3$ in \cite{end1}. In particular, they do not impose \emph{explicitily} any restriction on the rate of convergence of $(p_i,t_i)$ to $(p,T)$. Their analysis shows that, under the global Type~I assumption \eqref{eq.typeIupper}, this set coincides with the entire singular set $\Sigma$ and with the set $\Sigma_{\Rm}$ of points $p$ for which $\abs{\Rm(p,t)}_{g(t)}$ blows up at a Type~I rate (Theorem $1.2$ in \cite{end1}), meaning that one could ask for the existence of an essential blow-up sequence with $p_i=p$ for all $i\in\mathbb{N}$ in this case. Moreover, the analysis of the precise asymptotic behaviour of the neckpinch singularity developed by Angenent and Knopf \cite{ang} shows that, given a singular point $p$, the curvature tensor actually \emph{cannot} blow up along a sequence of regular space-time points $(p_i,t_i)$ with $(T-t_i)=\mathfrak{o}(d^2_{g(t_i)}(p_i,p))$. Such a sequence is in fact sent to infinity by rescaling parabolically with $(T-t_i)$. This suggests that the rate of convergence of any essential blow-up sequence should be such that \eqref{eq.convergenceofseq} holds.

The situation becomes more convoluted when the flow is not Type~I, and one needs to carefully check the balance between the different rates involved. A good example to have in mind is one of the degenerate neck-pinch solutions constructed by Angenent, Isenberg, and Knopf in \cite{ang2}. They prove that any blow-up sequence $(p,t_i)$ based at a fixed point $p$ in the smallest neck, and with rescaling factor $(T-t_i)$, converges to a shrinking cylinder while any blow-up sequence $(q,t_i)$ based at the tip $q$, with scaling factor given by the curvature at the tip at time $t_i$, converges to a Bryant soliton. From the first convergence result it is clear that $(T-t_i)=\mathfrak{o}(d^2_{g(t_i)}(p,q))$, which geometrically corresponds to the tip being sent to infinity for this blow-up sequence. Heuristically, we may think that a sequence of points converging to the tip will have curvature going to infinity at a Type~I or Type~II rate depending on how fast it converges to the tip. These argument suggests that suitable definitions of the Type~I and Type~II singular sets for general flows should involve the rate of convergence of the essential blow-up sequences, and the condition \eqref{eq.convergenceofseq} should give the right scale.

Before starting the proofs of the theorems from the introduction, we need the following adaptation of the ``shrinking and expanding balls lemmas'' in \cite{sim2,sim3} to our Riemann scale.

\begin{lemma}\label{shrinkingballslemma}
Suppose $(M,g(t))$ is a complete $n$-dimensional Ricci flow defined on $[0,T)$. Then there exists a constant $C_1=C_1(n)$ such that for any point $(p,t)\in M \times [0,T)$ we have the inclusion
\begin{equation}
B_{g(t)}(p,r_{\Rm}(p,t)) \supseteq B_{g(s)}(p,r_{\Rm}(p,t)-C_1^2 r_{\Rm}^{-1}(p,t)\abs{s-t}),
\end{equation}
for all $s\in(t-r_{\Rm}^2(p,t), t+r_{\Rm}^2(p,t))$.
\end{lemma}
\begin{proof}
We set $C_1=2\sqrt[4]{2/3} \sqrt{n-1}$ where $n$ is the dimension of the manifold $M$. Since $\Ric_{g(t)}\leq (n-1)r_{\Rm}^{-2}(p,t)g(t)$ on $\mathcal{P}(p,t,r_{\Rm})$, and the function $r_{\Rm}^{-2}(p,t)$ is constant, hence in particular continuous and integrable on the interval $[t-r_{\Rm}^2(p,t), t+r_{\Rm}^2(p,t)]$, we can appeal to Lemma 3.2 in \cite{sim2} (note that their constant $\beta/2$ is equal to our $C_1^2$) for any $s \in [t,t+r_{\Rm}^2(p,t))$ to obtain
\begin{equation*}
 B_{g(t)}(p,r_{\Rm}(p,t)) \supseteq B_{g(s)}(p,r_{\Rm}(p,t)-C_1^2 r_{\Rm}^{-1}(p,t)(s-t)).
\end{equation*}
In order to prove the inclusion for times $s \in (t-r_{\Rm}^2(p,t),t]$, we notice that due to the lower bound $\Ric_{g(t)}\geq - (n-1)r_{\Rm}^{-2}(p,t_0)g(t)$ on $\mathcal{P}(p,t,r_{\Rm})$, we can use Lemma 2.1 from \cite{sim3} to obtain
\begin{equation*}
 B_{g(t)}(p,r_{\Rm}(p,t)) \supseteq B_{g(s)}(p,r_{\Rm}(p,t)e^{-(n-1)r_{\Rm}^{-2}(p,t)(t-s)}).
\end{equation*}
Using the elementary inequalities $e^{-x}\geq 1-x$ and $C_1^2 \geq (n-1)$ we get
\begin{align*}
 B_{g(t)}(p,r_{\Rm}(p,t)) &\supseteq B_{g(s)}(p,r_{\Rm}(p,t)-(n-1)r_{\Rm}^{-1}(p,t)(t-s))\\
 &\supseteq B_{g(s)}(p,r_{\Rm}(p,t)- C_1^2 r_{\Rm}^{-1}(p,t)(t-s)).\qedhere
 \end{align*}
\end{proof}

We can now deduce the following estimate for the Riemann scale.

\begin{theorem}[Lipschitz-H\"older Continuity of Riemann Scale]\label{thm.lipschitzholder}
Suppose $(M,g(t))$ is a complete $n$-dimensional Ricci flow defined on $[0,T)$. Then for any pair of space-time points $(p,t)$ and $(q,s)$ we have
\begin{equation}\label{eq.lipschitzholder}
\abs{r_{\Rm}(p,t)-r_{\Rm}(q,s)}\leq \min \set{d_{g(t)}(p,q),d_{g(s)}(p,q)}+ C_1\abs{t-s}^{\frac{1}{2}},
\end{equation}
where $C_1$ is a constant depending only on the dimension $n$.
\end{theorem}

\begin{proof}
We first prove that the Riemann scale is Lipschitz continuous in space and then prove the H\"older continuity in time.

\emph{Step 1:} We show that for fixed $t\in[0,T)$ the Riemann scale $r_{\Rm}(\cdot,t)$ is $1$-Lipschitz continuous with respect to the metric $g(t)$. The proof of this part is a close adaptation of the analogous result in \cite{bamz}, but since its proof is as short as enlightening, we quickly recall it here. 

Suppose towards a contradiction that for some $t$ the function $r_{\Rm}(\cdot,t)$ is not $1$-Lipschitz with respect to $g(t)$, that is we can find $p$ and $q$ such that
\begin{equation*}
 r_{\Rm}(p,t)-r_{\Rm}(q,t)>d_{g(t)}(p,q).
\end{equation*}
Here we assumed without loss of generality that $r_{\Rm}(p,t) \geq r_{\Rm}(q,t)$. We then define $r := r_{\Rm}(p,t)-d_{g(t)}(p,q) > r_{\Rm}(q,t)$, so that $B_{g(t)}(q,r)\subseteq B_{g(t)}(p,r_{\Rm}(p,t))$ by the triangle inequality as well as $[t-r^2,t+r^2]\subseteq [t-r^2_{\Rm}(p,t),t+r^2_{\Rm}(p,t)]$, hence 
\begin{equation*}
\mathcal{P}(q,t,r) \subseteq \mathcal{P}(p,t,r_{\Rm}(p,t)).
\end{equation*} 
In particular, by definition of the Riemann scale at $(p,t)$, we have $\abs{\Rm}\leq r_{\Rm}^{-2}(p,t) \leq r^{-2}$ on $\mathcal{P}(q,t,r)$, and therefore $r_{\Rm}(q,t)\geq r$ by definition of the Riemann scale at $(q,t)$. This yields the desired contradiction with the definition of $r$.

\emph{Step 2:} We show that $r_{\Rm}(p,\cdot)$ is $\frac{1}{2}$-H\"older continuous with constant $C_1=C_1(n)$ being given by Lemma \ref{shrinkingballslemma}.

Arguing by contradiction, let us assume that for some fixed $p\in M$ there exist two times $t$ and $s$ in $[0,T)$ such that (assuming without loss of generality that $r_{\Rm}(p,t)\geq r_{\Rm}(p,s)$)
\begin{equation}\label{eq.nonHolder}
r_{\Rm}(p,t)-r_{\Rm}(p,s)>C_1\sqrt{\abs{t-s}},
\end{equation}
with $C_1=2\sqrt[4]{2/3} \sqrt{n-1}$ the constant from Lemma \ref{shrinkingballslemma}. Set $r := r_{\Rm}(p,t)-C_1\sqrt{\abs{t-s}} > r_{\Rm}(p,s)$. We claim that 
\begin{equation}\label{eq.Pinclusion}
\mathcal{P}(p,s,r) \subseteq \mathcal{P}(p,t,r_{\Rm}(p,t)).
\end{equation}
As above, if this claim is true, then we can deduce from the definition of $r_{\Rm}(p,t)$ that $\abs{\Rm}< r^{-2}_{\Rm}(p,t)<r^{-2}$
on $\mathcal{P}(p,s,r)$, therefore by definition of $r_{\Rm}(p,s)$ we obtain $r_{\Rm}(p,s)\geq r$, in contradiction with the definition of $r$, concluding the proof.

It remains to verify the claim \eqref{eq.Pinclusion}. We first check the time intervals. Note that we can consider the time intervals in $\R$ rather than $[0,T)$, therefore omitting the truncation at $0$ and $T$ as this inclusion clearly implies the one between the truncated intervals. By \eqref{eq.nonHolder}, and using $C_1\geq 1$, we can estimate
\begin{equation*}
2C_1 r_{\Rm}(p,t)\sqrt{\abs{t-s}} > 2C_1^2 \abs{t-s}\geq C_1^2 \abs{t-s}-t+s,
\end{equation*}
as well as
\begin{equation*}
2C_1 r_{\Rm}(p,t)\sqrt{\abs{t-s}} > 2C_1^2 \abs{t-s}\geq C_1^2 \abs{t-s}+t-s.
\end{equation*}
Therefore, we obtain
\begin{align*}
s+r^2=s+r_{\Rm}^2(p,t)+C_1^2\abs{t-s} -2C_1r_{\Rm}(p,t)\sqrt{\abs{t-s}}&\leq t+r_{\Rm}^2(p,t),\\
s-r^2=s-r_{\Rm}^2(p,t)-C_1^2\abs{t-s} +2C_1r_{\Rm}(p,t)\sqrt{\abs{t-s}}&\geq t-r_{\Rm}^2(p,t),
\end{align*}
and hence $(s-r^2,s+r^2) \subseteq (t-r_{\Rm}^2(p,t), t+r_{\Rm}^2(p,t))$. In order to prove the inclusion $B_{g(s)}(p,r) \subseteq B_{g(t)}(p,r_{\Rm}(p,t))$, we recall that Lemma \ref{shrinkingballslemma} implies
\begin{equation*}
 B_{g(t)}(p,r_{\Rm}(p,t)) \supseteq B_{g(s)}(p,r_{\Rm}(p,t)-C_1^2 r_{\Rm}^{-1}(p,t)\abs{t-s})
\end{equation*}
and hence we are done if we can show that $B_{g(s)}(p,r) \subseteq B_{g(s)}(p,r_{\Rm}(p,t)-C_1^2 r_{\Rm}^{-1}(p,t)\abs{t-s})$. To verify this, note that, again using \eqref{eq.nonHolder}, we can estimate
\begin{equation*}
C_1^2 r_{\Rm}^{-1}(p,t)\abs{t-s} < C_1\sqrt{\abs{t-s}}
\end{equation*}
and therefore
\begin{equation*}
r_{\Rm}(p,t)-C_1^2 r_{\Rm}^{-1}(p,t)\abs{t-s} > r_{\Rm}(p,t)-C_1\sqrt{\abs{t-s}}=r.
\end{equation*}
This finishes the proof of Step 2.

\emph{Step 3:} Finishing the proof of Theorem \ref{thm.lipschitzholder} is now straight-forward: From the above two steps, an easy application of the triangle inequality both in space and time yields \eqref{eq.lipschitzholder} for $C_1=2\sqrt[4]{2/3} \sqrt{n-1}$.
\end{proof}

A first corollary of this theorem gives upper and lower bounds for $r^{-2}_{\Rm}$ on a parabolic cylinder in terms of its value at the center of this cylinder.
\begin{corollary}[Local Harnack-Type~Inequality]\label{harnackcor}
Suppose $(M,g(t))$ is a complete $n$-dimensional Ricci flow defined on $[0,T)$, and let $(p,t)\in M \times (0,T)$. Then there exists a constant $a_1=a_1(n)\in(0,1)$ such that
\begin{equation}\label{eq.Harnack}
\frac{1}{4} r^{-2}_{\Rm}(p,t)\leq r^{-2}_{\Rm}(\cdot,\cdot) \leq 4 r^{-2}_{\Rm}(p,t) \quad \text{on } \mathcal{P}(p,t,a_1\,r_{\Rm}(p,t))
\end{equation}
\end{corollary}

\begin{proof}
Let $(q,s)\in \mathcal{P}(p,t,a_1\,r_{\Rm})$, where $a_1<1$ will be determined later. We deduce from Theorem \ref{thm.lipschitzholder} that
\begin{equation*}
\abs{r_{\Rm}(p,t)-r_{\Rm}(q,s)}\leq d_{g(t)}(p,q)+ C_1\abs{t-s}^{\frac{1}{2}} \leq a_1\,r_{\Rm}(p,t) +C_1a_1\,r_{\Rm}(p,t).
\end{equation*}
Rearranging this inequality as
\begin{equation*}
 (1+a_1+a_1 C_1) r_{\Rm}(p,t) \geq r_{\Rm}(q,s)
\end{equation*}
we see that the first inequality in \eqref{eq.Harnack} is implied by
\begin{equation*}
 1+a_1 +a_1 C_1 \leq 2 \iff a_1 \leq \frac{1}{1+C_1}.
\end{equation*}
Similarly, rearranging differently, we see that
\begin{equation*}
 (1-a_1-a_1 C_1) r_{\Rm}(p,t) \leq r_{\Rm}(q,s),
\end{equation*}
therefore for the second inequality in \eqref{eq.Harnack} it is sufficient to impose
\begin{equation*}
 1-a_1-a_1 C_1 \geq \frac{1}{2} \iff a_1 \leq \frac{1}{2(1+C_1)}.
\end{equation*}
The claim hence follows by setting $a_1:= \frac{1}{2(1+C_1)}$ with $C_1$ as above.
\end{proof}

We are now ready to prove Theorem \ref{thm.singRiemannscale}. Maybe slightly paradoxically, the most subtle argument is actually needed in the proof of Part $i)$, because the definition of singular and regular points is not set in a parabolic way. Proving the other two statements instead is easier as our definitions of Type~I and Type~II singular points are already in terms of parabolic neighbourhoods. 

\begin{proof}[Proof of Theorem \ref{thm.singRiemannscale}]
Let $(M,g(t))$ be a Ricci flow defined on $[0,T)$, $T<\infty$ satisfying \eqref{eq.bddcurv} (possibly incomplete) and let $p\in M$.
\begin{enumerate}[i)]
\item We prove the following equivalent statement to the claim in the theorem: $p$ is singular if and only if $\liminf_{t\nearrow T} r_{\Rm}(p,t)=0$.

To this end, suppose first that $p$ is regular, and let $U$ be a neighbourhood of $p$ and $C>0$ a constant such that $\abs{\Rm}\leq C$ on $U \times [0,T)$. Let $r_0>0$ be such that $B_{g(0)}(p,r_0)\subset \subset U$. The standard multiplicative distance distortion estimate gives the containment
\begin{equation}
B_{g(0)}(p,r_0) \supseteq B_{g(t)}(p,r_0e^{-(n-1)Ct}) \supseteq B_{g(t)}(p,r_1),
\end{equation}
where $r_1 := r_0e^{-(n-1)CT}>0$. From this we deduce $r_{\Rm}(p,t) \geq \min \{ C^{-1/2},r_1 \}>0$.

For the converse statement, suppose that there exists some constant $\delta>0$ such that $\liminf_{t\nearrow T} r_{\Rm}(p,t) >\delta$. By the curvature boundedness assumption, we can choose this $\delta$ so that we have $r_{\Rm}(p,t) >\delta$ for every $t \in [0,T)$; by definition this means that for every $t \in [\delta^2,T)$ we have $\abs{\Rm} \leq \delta^{-2}$ on $\mathcal{P}(p,t,r_{\Rm}) \supseteq \mathcal{P}(p,t,\delta)$. Fixing any $t$ in this range, we claim that there exists a constant $r_2=r_2(\delta,n,T)$ such that 
\begin{equation}\label{eq.r2}
B_{g(t)}(p,\delta) \supseteq B_{g(\delta^2)}(p,r_2).
\end{equation} 
In particular, if this claim holds, we have found a (fixed) neighbourhood of $p$ on which the curvature remains bounded (by $\delta^{-2}$), and we can conclude the proof. Hence it remains to prove the claim. In order to do so, for the $t$ we fixed above, we set $k_0=k_0(t)$ to be the smallest integer such that $t-(k_0+1)\delta^2 <0$. Notice that $k_0 \leq T/\delta^2$. We are going to implement an iterative scheme of inclusions using the distance distortion estimates in any of the cylinders considered. From the bound on the curvature we obtain
\begin{align*}
B_{g(t)}(p,\delta) &\supseteq B_{g(t-\delta^2)}\big(p,\delta e^{-(n-1)\delta^{-2} \cdot \delta^2}\big)=B_{g(t-\delta^2)}\big(p,\delta e^{-(n-1)}\big)\\
&\supseteq B_{g(t-2\delta^2)}\big(p,\delta (e^{-(n-1)})^2\big)=B_{g(t-2\delta^2)}\big(p,\delta e^{-2(n-1)}\big)\\
&\supseteq .... \supseteq B_{g(t-(k_0-1)\delta^2)}\big(p,\delta e^{-(k_0-1)(n-1)}\big)\\
&\supseteq B_{g(\delta^2)}\big(p,\delta e^{-(n-1)[(k_0-1)+\delta^{-2}(t-(k_0-1)\delta^2-\delta^2)]}\big)\\
&\supseteq B_{g(\delta^2)}\big(p,\delta e^{-k_0(n-1)}\big) \supseteq B_{g(\delta^2)}\big(p,\delta e^{-(n-1)(T/\delta^2)}\big),
\end{align*}
and \eqref{eq.r2} follows by defining $r_2 := \delta e^{-(n-1)(T/\delta^2)}$.

\item Let us first assume the lower bound in (\ref{eq.typeI}) and let $t_i\nearrow T$ be a sequence realising the $\limsup$. Then, setting $r := \max \set{r_I, 1/\sqrt{c_I}}$, we have for sufficiently large $i$
\begin{equation*}
\frac{1}{r^2(T-t_i)} \leq \frac{c_I}{T-t_i} < \sup_{\mathcal{P}(p,t_i,r_I\sqrt{T-t_i})} \abs{\Rm} \leq \sup_{\mathcal{P}(p,t_i,r\sqrt{T-t_i})} \abs{\Rm}.
\end{equation*}
This means that we must have $r_{\Rm}(p,t_i) < r\sqrt{T-t_i}$ for all sufficiently large $i$, because if it was $r_{\Rm}(p,t_j) \geq r\sqrt{T-t_j} =: r'$ for some $j$ large enough for the above inequality to hold, we would obtain
\begin{equation*}
\frac{1}{r^2(T-t_j)} < \sup_{\mathcal{P}(p,t_j,r')} \abs{\Rm} \leq \frac{1}{(r')^2}=\frac{1}{r^2(T-t_j)},
\end{equation*}
a contradiction. We therefore conclude that 
\begin{equation*}
r_{\Rm}(p,t_i)/\sqrt{T-t_i} < \max \set{r_I,1/\sqrt{c_I}},
\end{equation*} 
or equivalently
\begin{equation}\label{eq.alowerbound}
(T-t_i)r^{-2}_{\Rm}(p,t_i)> \min \set{r_I^{-2},c_I} =: \widetilde{c}_I>0,
\end{equation} 
which gives the lower bound we claimed. Supposing instead the upper bound of (\ref{eq.typeI}), for $\e>0$, let us set $r := \min \set{r_I,1/\sqrt{C_I+\e}}$, so that we can compute for any $t$ sufficiently close to $T$ that
\begin{equation*}
\frac{1}{r^2(T-t)} \geq \frac{C_I+\e}{T-t} \geq \sup_{\mathcal{P}(p,t,r_I\sqrt{T-t})} \abs{\Rm} \geq \sup_{\mathcal{P}(p,t,r\sqrt{T-t})} \abs{\Rm},
\end{equation*}
so by definition we have $r_{\Rm}(p,t) \geq r\sqrt{T-t}$. In particular, letting $\e\searrow 0$, we conclude that 
\begin{equation*}
\liminf_{t \nearrow T} r_{\Rm}(p,t)/\sqrt{T-t} \geq \min \set{r_I, 1/\sqrt{C_I}},
\end{equation*}
or equivalently
\begin{equation}\label{eq.anupperbound}
\limsup_{t \nearrow T} (T-t)r^{-2}_{\Rm}(p,t) \leq \max \set{r_I^{-2}, C_I} =: \widetilde{C}_I.
\end{equation}

To prove the converse statement, assume now that $\widetilde{c}_I < \limsup_{t\nearrow T}\, (T-t) r_{\Rm}^{-2}(p,t)$ for some constant $0<\widetilde{c}_I$. Set $r_I:= (\widetilde{c}_I)^{-\frac{1}{2}}$ and $c_I:=\widetilde{c}_I$. We claim that for this choice, the lower bound in \eqref{eq.typeI} must hold. If not, then we have
\begin{equation*}
 \limsup_{t\nearrow T} \sup_{B_{\widetilde{g}_t(0)}(p,r_I)\times(-r_I^2,r_I^2)} \,\abs{\Rm_{\widetilde{g}_t}}_{\widetilde{g}_t} \leq c_I,
\end{equation*}
that is, we have an upper bound as in \eqref{eq.typeI}, but with $C_I$ now replaced by $c_I$. By what we have just proved above, we thus get the analogue of \eqref{eq.anupperbound}, namely
\begin{equation*}
\limsup_{t \nearrow T} (T-t)r^{-2}_{\Rm}(p,t) \leq \max \set{r_I^{-2}, c_I} = \widetilde{c}_I,
\end{equation*}
which is the desired contradiction. Finally, for the last remaining statement, assume that $\limsup_{t\nearrow T}\, (T-t) r_{\Rm}^{-2}(p,t) \leq \widetilde{C}_I$. Then given any $\e$, we have for every $t$ close enough to $T$ that 
\begin{equation*}
r_{\Rm}(p,t) > \sqrt{\frac{T-t}{\widetilde{C}_I+\e}} =:r.
\end{equation*} 
The definition of the Riemann scale thus implies
\begin{equation}
 \sup_{\mathcal{P}(p,t,r)} \abs{\Rm} \leq \sup_{\mathcal{P}(p,t,r_{\Rm})} \abs{\Rm}= r^{-2}_{\Rm}(p,t) \leq r^{-2}=\frac{\widetilde{C}_I+\e}{T-t},
\end{equation}
which for $\e\searrow 0$ and after rescaling parabolically gives the upper bound in \eqref{eq.typeI} with $r_I:= (\widetilde{C}_I)^{-\frac{1}{2}}$ and $C_I:=\widetilde{C}_I$.

\item The argument used in the first paragraph of the proof of Part $ii)$ also works with Type~II points, that is verifying (\ref{eq.typeII}), and one shows that for every $r>0$ and every $c_I>0$ the analogue of \eqref{eq.alowerbound} holds, namely there exists a sequence of $t_i\nearrow T$ such that
\begin{equation*}
(T-t_i)r^{-2}_{\Rm}(p,t_i)> \min \set{r^{-2},c_I}.
\end{equation*} 
We can therefore let $r\searrow 0$ and $c_I\nearrow\infty$ to conclude that $\limsup_{t\nearrow T}\, (T-t) r_{\Rm}^{-2}(p,t)=\infty$.

Conversely, suppose we have $\limsup_{t\nearrow T}\, (T-t) r_{\Rm}^{-2}(p,t)=\infty$. If \eqref{eq.typeII} does not hold for all $r>0$, then there exists some $r_I>0$ and some constant $C_I$ such that the upper bound in \eqref{eq.typeI} holds. But we have just proven that this implies $\limsup_{t\nearrow T}\, (T-t) r_{\Rm}^{-2}(p,t)\leq \widetilde{C}_I$, a contradiction. This finishes the proof of the theorem.\qedhere
\end{enumerate}
\end{proof}

With a slight modification of Part $i)$ of Theorem \ref{thm.singRiemannscale} we obtain a non-oscillation result which in particular, says that, to a certain extent, the curvature of a Ricci flow cannot oscillate between a Type~I rate and a lower rate arbitrarily close to the singular time.

\begin{corollary}[Type~I non-Oscillation]\label{typeIoscicor}
Suppose $(M,g(t))$ is a Ricci flow defined on a finite time interval $[0,T)$, satisfying either the bounded curvature condition \eqref{eq.bddcurv} or having complete time slices for all $t \in [0,T)$. Then $p \in \Sigma$ if and only if
\begin{equation}\label{eq.typeIoscicor}
r^{-2}_{\Rm}(p,t) > \frac{1}{T-t}, \quad \forall t\in[0,T).
\end{equation}
\end{corollary}

\begin{proof}
If \eqref{eq.typeIoscicor} holds, then by Part $i)$ of Theorem \ref{thm.singRiemannscale} the point $p$ must be singular. Conversely, assume towards a contradiction that $p\in\Sigma$, but that
\begin{equation*}
\delta_0 := r_{\Rm}(p,t_0) \geq \sqrt{T-t_0}, \quad \text{for some }t_0\in[0,T).
\end{equation*}
By definition of the Riemann scale, this means in particular that
\begin{equation}\label{eq.boundfromdef}
\abs{\Rm} \leq \delta_0^{-2}, \quad \text{on }\mathcal{P}(p,t_0,\delta_0) \supseteq B_{g(t_0)}(p,\delta_0)\times [t_0,T).
\end{equation}
Set $U:=B_{g(t_0)}(p,\delta_0)$. If the flow has bounded curvature, we can choose $\delta\leq \delta_0$ such that $\abs{\Rm} \leq \delta^{-2}$ on $U\times[0,t_0]$. Otherwise, by the completeness assumption, the continuous function $\abs{\Rm}$ is less than or equal to $\delta^{-2}$ for some $\delta \leq \delta_0$ on the compact set $U\times [0,t_0]$. Hence in any case, combining this with \eqref{eq.boundfromdef}, $\abs{\Rm} \leq \delta^{-2}$ on $U\times[0,T)$ and therefore by definition $p\in \mathfrak{Reg}$, yielding the desired contradiction.
\end{proof}

The decomposition of the singular set given by Theorem \ref{thm.singdecomposition} follows very easily now.
\begin{proof}[Proof of Theorem \ref{thm.singdecomposition}]
From the equivalent definitions of Type~I and Type~II singular points given by Theorem \ref{thm.singRiemannscale}, it is clear that if $p\in\Sigma_I\cup\Sigma_{II}$ then $p\in\Sigma$.

Conversely, if $p\in\Sigma$, then by Corollary \ref{typeIoscicor} we have
\begin{equation*}
r^{-2}_{\Rm}(p,t) > \frac{1}{T-t}, \quad \forall t\in[0,T),
\end{equation*}
and therefore in particular
\begin{equation*}
\limsup_{t\nearrow T}\, (T-t) r_{\Rm}^{-2}(p,t) \geq 1> 0.
\end{equation*}
By Theorem \ref{thm.singRiemannscale}, we therefore have $p\in\Sigma_I\cup\Sigma_{II}$.
\end{proof}

As an immediate consequence of Theorem \ref{thm.singdecomposition}, we find the following corollary which is equivalent to Theorem $3.2$ in \cite{end1}.
\begin{corollary}
Let $(M,g(t))$ be a Type~I Ricci flow on $[0,T)$, $T<+\infty$. Then $\Sigma=\Sigma_I$.
\end{corollary} 

Next, we give some evidence for the following conjecture. 

\begin{conjecture}\label{con1}
In Theorem \ref{thm.singRiemannscale}, we can replace each instance of $\limsup_{t\nearrow T}$ with $\liminf_{t\nearrow T}$.
\end{conjecture}

A full proof of this conjecture would require ruling out significant oscillations of the curvature between Type~I and Type~II rates. Intuitively, such an oscillatory behaviour would be in extreme contrast with the parabolic nature of the Ricci flow; nevertheless, it is not clear to the authors how to prevent it using only the differential inequalities on the curvature tensor given by Shi's estimates and even the H\"older continuity proved for the Riemannian scale is not strong enough to rule this out completely. The following result is the best we can currently prove in this direction. It studies the convergence rates of an essential blow-up sequence $(p_i,t_i)$ along which the curvature blows up at a Type~I rate but which converges to a Type~II singular point.
\begin{proposition}\label{typeIIosciprop}
Let $(M,g(t))$ be a complete Ricci flow on $[0,T)$, $T<\infty$, satisfying \eqref{eq.curvatureblowup} and \eqref{eq.bddcurv}. Suppose $p\in \Sigma_{II}$, and let $(p_i,t_i)$ be an essential blow-up sequence, with $p_i \rightarrow p$ in the topology of $(M,g(0))$, $t_i\nearrow T$, and such that $r_{\Rm}^{-2}$ blows up at a Type~I rate along $(p_i,t_i)$, i.e.
\begin{equation*}
r_{\Rm}^{-2}(p_i,t_i) \leq \frac{1}{m^2 (T-t_i)}
\end{equation*}
for some $m\in (0,1)$. Then there exists $\delta=\delta(n,m)>0$ such that
\begin{equation*}
 \limsup_{i \rightarrow \infty} \frac{T-t_{i-1}}{T-t_i} \geq 1+\delta \qquad \text{or} \qquad \limsup_{i \rightarrow +\infty} \frac{d^2_{g(t_i)}(p,p_i)}{T-t_i}\geq \delta.
\end{equation*}
\end{proposition}

\begin{proof}
Suppose instead that
\begin{equation}\label{lineartimesequence}
\limsup_{i \rightarrow +\infty} \frac{T-t_{i-1}}{T-t_i}=1 \qquad \text{and} \qquad \limsup_{i \rightarrow +\infty} \frac{d^2_{g(t_i)}(p,p_i)}{T-t_i}= 0.
\end{equation}
Since by hypothesis $r_{\Rm}^2(p_i,t_i) \geq m^2 (T-t_i)$, we infer
\begin{equation*}
\bigcup_{i \in \mathbb{N}} \big(t_i-a_1\,r_{\Rm}^2(p_i,t_i),t_i+a_1\,r_{\Rm}^2(p_i,t_i)\big) \supseteq \bigcup_{i \in \mathbb{N}} \big(t_i-m^2 a_1^2(T-t_i),t_i+m^2 a_1^2(T-t_i)\big).
\end{equation*}
We first check that the set on the right hand side contains $(t_{i_1},T)$ for some $i_1$. Indeed, the intervals in consideration overlap definitively as
\begin{equation*}
t_{i-1} +m^2 a_1^2(T-t_{i-1})\geq t_i-m^2 a_1^2(T-t_i) \iff T-t_{i-1} \leq \frac{1+m^2 a_1^2}{1-m^2 a_1^2} (T-t_i),
\end{equation*}
which is satisfied by (\ref{lineartimesequence}) for any $i$ large enough.
On the other hand, we also have
\begin{equation*}
 B_{g(t_i)}(p_i,a_1\,r_{\Rm}(p_i,t_i))\supseteq B_{g(t_i)}(p_i,a_1m\sqrt{T-t_i}) \ni p
\end{equation*}
for $i$ large enough by (\ref{lineartimesequence}). Therefore, for any sequence $\xi_j \nearrow T$ and for every $j$ large enough, there exists $i=i(j)$ such that
\begin{equation*}
 (p,\xi_j)\in \mathcal{P}(p_i,t_i,a_1\,r_{\Rm}(p_i,t_i)).
\end{equation*}
We can therefore appeal to the local Harnack inequality of Corollary \ref{harnackcor} to compare the values $r_{\Rm}^{2}(p,\xi_j)$ to $r_{\Rm}^{2}(p_i,t_i)$, where $i$ depends on $j$.

Since $p \in \Sigma_{II}$, by Theorem \ref{thm.singRiemannscale} we can pick a sequence of times $\xi_j\nearrow T$ such that $r^2_{\Rm}(p,\xi_j)=\mathfrak{o}(T-\xi_j)$. From the discussion above we can apply Corollary \ref{harnackcor} to obtain
\begin{equation}
\mathfrak{o}(T-\xi_j)=4 r_{\Rm}^2(p,\xi_j) \geq r_{\Rm}^2(p_i,t_{i(j)}) \geq m^2 (T-t_{i(j)}) \geq \frac{m^2}{1+m^2 a_1^2} (T-\xi_j),
\end{equation}
which gives a contradiction for $j$ large enough.
\end{proof}
\begin{remark}
The proposition above ensures in particular that for a point $p \in \Sigma_{II}$, the existence of a sequence $t_i\nearrow T$ along which $r_{\Rm}^{-2}(p,t_i)$ blows up at a Type~I rate forces the convergence to be exponential. We think that this property is in contrast to the Lipschitz continuity of $r_{\Rm}(p,\cdot)$.  It could be instructive to notice that the same is a priori not true if we swap the roles of Type~I and Type~II. Indeed, consider the function
\begin{equation}
f(t) = (T-t)^2+(T-t)\abs{\sin(\ln(T-t))}.
\end{equation}
This is a Lipschitz function, such that $f(t)=(T-t_k)^2$ on $t_k=T-\exp (-k \pi)$, and $f(t)\sim T-t$ elsewhere as $t \sim T$.

From a more optimistic point of view, if Conjecture \ref{con1} holds, it is natural to ask whether we must have 
\begin{equation*}
\liminf_{t \nearrow T}r^{-2}_{\Rm}(p,t)(T-t)=\limsup_{t \nearrow T}r^{-2}_{\Rm}(p,t)(T-t),
\end{equation*}
at any singular point $p \in \Sigma$. This would be coherent with the examples in the literature.
\end{remark}

For completeness, let us also study the convergence of an essential blow-up sequences along which the curvature blows up at a Type~II rate but which converges to a Type~I singular point $p$. From Definition \ref{def.singpoints}, it is clear that this convergence cannot be too fast. The following proposition makes this more precise, giving an explicit relation between the blow-up rate of $r^{-2}_{\Rm}(p,\cdot)$ and the convergence rate.
\begin{proposition}
Let $(M,g(t))$ be a complete Ricci flow on $[0,T)$, $T<\infty$, satisfying \eqref{eq.curvatureblowup} and \eqref{eq.bddcurv}. Let $p\in\Sigma_I$, meaning in particular that there exists $m\in(0,1)$ such that $r_{\Rm}^{-2}(p,t) \leq \frac{1}{m^2 (T-t)}$ for all $t\in[0,T)$. Suppose further that $(p_i,t_i)$ is an essential blow-up sequence, with $p_i \rightarrow p$ in the topology of $(M,g(0))$, $t_i\nearrow T$, and such that $r_{\Rm}^{-2}$ blows up at a Type~II rate along $(p_i,t_i)$, i.e.
\begin{equation*}
r^2_{\Rm}(p_i,t_i)=\mathfrak{o}(T-t_i).
\end{equation*} 
Then we obtain
\begin{equation*}
\limsup_{i \rightarrow \infty} \frac{d^2_{g(t_i)}(p,p_i)}{T-t_i} \geq m^2 a_1^2,
\end{equation*}
where $a_1=a_1(n)$ is the constant from Corollary \ref{harnackcor}.
\end{proposition}

\begin{proof}
Suppose by contradiction that
\begin{equation*}
 \limsup_{i \rightarrow \infty} \frac{d^2_{g(t_i)}(p,p_i)}{T-t_i} <m^2 a_1^2,
\end{equation*}
so that we obtain
\begin{equation*}
 B_{g(t_i)}(p,a_1\,r_{\Rm}(p,t_i))\supseteq B_{g(t_i)}(p,a_1 m\sqrt{T-t_i}) \ni p_i
\end{equation*}
for $i$ large enough. This means in particular that
\begin{equation*}
 (p_i,t_i)\in \mathcal{P}(p,t_i,a_1\,r_{\Rm}(p,t_i)),
\end{equation*}
for $i$ large enough and we can therefore use the local Harnack-type inequality of Corollary \ref{harnackcor} to compare the values $r^2_{\Rm}(p_i,t_i)$ to $r^2_{\Rm}(p,t_i)$ to obtain
\begin{equation*}
\mathfrak{o}(T-t_i)=4 r^2_{\Rm}(p_i,t_i) \geq r^2_{\Rm}(p,t_i) \geq m^2 (T-t_i)
\end{equation*}
which gives a contradiction for $i$ large enough.
\end{proof}

We conclude this section by establishing a link between the two different Riemann scales defined in Definition \ref{def.riemannscale}. As already said, these two scales are equivalent for bounded scalar curvature Ricci flows in view of the pseudolocality result in Proposition $3.2$ of \cite{bam}. For a general Ricci flow, one should not expect such a strong result, but a weaker infinitesimal analogue holds true.

\begin{proposition}[Characterisation of Singular Set using Fixed Time Slice Scale]\label{prop.singRiemannscale}
Let $(M,g(t))$ be a complete Ricci flow on $[0,T)$, $T<\infty$, satisfying \eqref{eq.curvatureblowup} and \eqref{eq.bddcurv}. Then $p \in \Sigma$ if and only if $\liminf_{t \nearrow T}\widetilde{r}_{\Rm}(p,t)=0$.
\end{proposition}

\begin{proof}
The implication $p \in \mathfrak{Reg}\Rightarrow \liminf_{t \nearrow T}\widetilde{r}_{\Rm}(p,t)>0$ follows the exact same lines as in the proof of Part $i)$ of Theorem \ref{thm.singRiemannscale}.

Conversely, assume that for a point $p \in M$ there exists some constant $\delta>0$ such that $\liminf_{t\nearrow T} \widetilde{r}_{\Rm}(p,t) >\delta$. Since the flow has bounded curvature, for a possibly smaller $\delta$ we have $\widetilde{r}_{\Rm}(p,t) >\delta$ for every $t \in [0,T)$; by definition of time-slice Riemann scale, this means that for every $t \in [0,T)$ we have $\abs{\Rm} \leq \delta^{-2}$ on $B_{g(t)}(p,\widetilde{r}_{\Rm}(p,t)) \supseteq B_{g(t)}(p,\delta)$. Thus we obtain 
\begin{equation*}
\abs{\Rm} \leq \delta^{-2} \ \ \text{on} \ \bigcup_{t \in [0,T)} B_{g(t)}(p,\delta) \times \set{t}.
\end{equation*}
We claim that there exists a constant $a_4=a_4(n)$ such that for any time $t_0 \in [0,T)$ we have $r_{\Rm}(p,t_0) >a_4 \delta$. Once we have proven this, we infer that $\liminf_{t\nearrow T} r_{\Rm}(p,t) \ge a_4 \delta$, and we can conclude the proof thanks to Theorem \ref{thm.singRiemannscale}, Part $i)$. Using the Expanding balls Lemma $3.1$ in \cite{sim2} with $R=\delta$ and $r=\frac{\delta}{2}$, the lower bound $\Ric_{g(t)} \ge -(n-1)\delta^{-2} g(t)$ on the balls $B_{g(t)}(p,\delta)$ ensures
\begin{equation*}
B_{g(t_0)}\big(p,\tfrac{\delta}{2}\big) \subseteq B_{g(t)}(p,\delta), \quad \forall t \in \Big[ t_0,\min \Big\{ t_0+\tfrac{\delta^2}{n-1} \log(2),T \Big\} \Big).
\end{equation*}
On the other hand, the upper bound $\Ric_{g(t)} \le (n-1)\delta^{-2} g(t)$ on the balls $B_{g(t)}(p,\delta)$ guarantees that we are in the hypothesis of the Shrinking balls Lemma $3.2$ in \cite{sim2} with our $t$ being their initial time $0$, $r=\delta$ and $f=\delta^{-1}$, so that we obtain
\begin{equation*}
B_{g(t_0)}\big(p,\tfrac{\delta}{2}\big) \subseteq B_{g(t_0)}(p,\delta-C_1^2 (t_0-t) \delta^{-1}) \subseteq B_{g(t)}(p,\delta), \quad \forall t \in \Big( \max \Big\{ 0,t_0-\tfrac{\delta^2}{2 C_1^2} \Big\},t_0 \Big].
\end{equation*}
Recall that their constant $\beta/2$ is equal to our $C_1^2$. Therefore, we have obtained the inclusion
\begin{equation*}
\mathcal{P}(p,t_0,a_4 \delta) \subseteq \bigcup_{t} B_{g(t)}(p,\delta) \times \set{t},
\end{equation*}
where the union is taken over $t\in  (\max \set{t_0-a_4^2 \delta^2,0},\min \set{t_0+a_4^2 \delta^2,T})$ and where $a_4=a_4(n):=\min \big\{\sqrt{\frac{\log(2)}{n-1}},\frac{1}{\sqrt{2} C_1},\frac{1}{2}\big\}$. By definition of the Riemann scale, we see that $r_{\Rm}(p,t_0) > a_4 \delta$, as we wanted to prove.
\end{proof}


\section{Integral Characterisation of the Singular Sets}\label{sec.mixedintegral}

Let us start with recalling the fundamental Noncollapsing Theorem of Perelman as well as an extension of it (see \cite{ye2}). We need the following definition.

\begin{definition}[$\kappa$-Noncollapsing at Curvature Scales]\label{def.noncollapsed}
Let $\kappa>0$. We say that a Ricci flow $g(t)$ is $\kappa$-noncollapsed on the scale $\varrho$ if every metric ball $B_{g(t)}(p,r)$ of radius $r < \varrho$ that satisfies for every $(x,t)\in B_{g(t)}(p,r) \times (t-r^2,t]$ the curvature bound $\abs{\Rm(x,t)}_{g(t)} \leq r^{-2}$, has volume at least $\kappa r^n$. We say that the flow is $\kappa$-noncollapsed on the scale $\varrho$ relative to the scalar curvature if for every metric ball $B_{g(t)}(p,r)$ of radius $r < \varrho$ that satisfies the scalar curvature bound $\abs{\Sc(\cdot,t)}_{g(t)} \leq n(n-1)r^{-2}$ on $B_{g(t)}(p,r)$, has volume at least $\kappa r^n$.
\end{definition}

For a complete Ricci flow $(M,g(t))$ defined on a finite time interval $[0,T)$, satisfying \eqref{eq.bddcurv} and a lower injectivity radius bound at the initial time, Perelman's Noncollapsing Theorem \cite{per} guarantees for every $\varrho$ the existence of a constant $\kappa=\kappa(n,g(0),T,\varrho)$ such that the flow is $\kappa$-noncollapsed on the scale $\varrho$.
Furthermore, if the flows has uniformly bounded scalar curvature, $\abs{\Sc}\le n(n-1)R_0$ on $M \times [0,T)$, it is $\kappa_1$-noncollapsed on the scale $\varrho = \min\{R_0^{-1/2},\sqrt{T}\}$ relative to the scalar curvature. In fact, by Aubin's classical result we have bounds on the Sobolev constants of the initial metric $g(0)$ (see \cite{aub}). These bounds extend to later times thanks to Lemma A.3 in \cite{z07}, which yields the claimed statement using Lemma A.4 in the same paper. For another approach to these noncollapsing result see \cite{ye2}. In particular, for any $r_0$ smaller than $\varrho$, we have the volume bound
\begin{equation}
\mu_{g(t_0)} \Big( B_{g(t_0)}(p_0,r_0)\Big) \ge \kappa_1 r_0^n.
\end{equation}
Here $\kappa_1$ depends on the dimension $n$, the initial metric $g(0)$, the time $T$ and the scalar curvature bound $R_0$.

The Harnack-type inequality proved in Corollary \ref{harnackcor} implies an integral concentration of the curvature. 

\begin{theorem}[Integral Curvature Concentration]\label{riemannscaleintegralbounds}
Let $\kappa>0$ and let $(\alpha,\beta)$ be an optimal pair of integrability exponents in the sense of Definition \ref{def.optimal}. Then there exist constants $C_2=C_2(n,\kappa,\alpha)>0$ and $C_3=C_3(n)$ such that the following holds. Let $(M,g(t))$ be a complete Ricci flow defined on $[0,T)$, $T<\infty$, which is $\kappa$-non-local-collapsed on a scale $\varrho$ in the sense of Definition \ref{def.noncollapsed}. Then for a space-time point $(p,t) \in \Sigma \times (T-\varrho^2,T)$, we have the integral bounds
\begin{equation}\label{eq.intbounds}
C_2 \leq \norm{r_{\Rm}^{-2}}_{\alpha,\beta,\mathcal{P}(p,t,a_1\,r_{\Rm})} \leq C_3,
\end{equation}
where $a_1\in(0,1)$ is the constant from Corollary \ref{harnackcor}.
\end{theorem}

This can be seen as an $\e$-regularity theorem since the lower bound in \eqref{eq.intbounds} shows that if $\norm{r_{\Rm}^{-2}}_{\alpha,\beta,\mathcal{P}(p,t,a_1\,r_{\Rm})} \leq \e< C_2$ as $t\to T$, then $p$ must be a regular point. 
\begin{proof}
Since $p\in\Sigma$, Corollary \ref{typeIoscicor} implies that $r^2_{\Rm}(p,t) \leq (T-t)$. On the one hand this implies that $r_{\Rm}(p,t)<\varrho$ for large enough $t$ so that we can use the $\kappa$-noncollapsing property for balls of radius $r\leq r_{\Rm}(p,t)$. On the other hand, it also shows that 
\begin{equation}\label{eq.intinclusion}
\begin{aligned}
t+a_1^2 r^2_{\Rm}(p,t) \leq t+a_1^2 (T-t)<t+\frac{1}{2}(T-t)<T,\\
t-a_1^2 r^2_{\Rm}(p,t) \geq t-a_1^2 (T-t)>t-\frac{1}{2}(T-t)>0.
\end{aligned}
\end{equation}
In particular, we have 
\begin{equation}\label{eq.nontruncated}
\mathcal{P}(p,t,a_1\,r_{\Rm}(p,t)) = B_{g(t)}(p,a_1\,r_{\Rm}(p,t)) \times (t-a_1^2 r^2_{\Rm}(p,t),t+a_1^2 r^2_{\Rm}(p,t))
\end{equation}
and we do not need to worry about the truncation in \eqref{eq.paraboliccylinderdef}. 

We now first prove the upper bound. By definition of the Riemann scale, we have a Riemann upper bound on $\mathcal{P}(p,t,r_{\Rm})$, so after rescaling and using the Bishop-Gromov inequality we obtain
\begin{equation*}
\mu_{g(t)}(B_{g(t)}(p,a_1\,r_{\Rm}(p,t)))\leq a_1^n r_{\Rm}^n(p,t) \mu_{g_{hyp}}(B_{hyp})=a_1^n C_H(n) r_{\Rm}^n(p,t),
\end{equation*}
where $B_{hyp}$ denotes a unitary ball in the hyperbolic space, and $C_H(n)=\int_0^1{\sinh^{n-1}(s) ds}$ is its volume.
Considering the evolution equation of the volume element under Ricci flow, and using the curvature bound in the region $\mathcal{P}(p,t,r_{\Rm})$, we deduce that for every time $s \in (t-a_1^2 r^2_{\Rm}(p,t),t+a_1^2 r^2_{\Rm}(p,t))$
\begin{equation}\label{eq.volumeupper}
\begin{aligned}
\mu_{g(s)}(B_{g(t)}(p,a_1\,r_{\Rm}(p,t))) &\leq e^{n(n-1)r_{\Rm}^{-2}(p,t)(t-s)} \mu_{g(t)}(B_{g(t)}(p,a_1\,r_{\Rm}(p,t)))\\
&\leq  e^{n(n-1)a_1^2} a_1^n C_H(n) r_{\Rm}^n(p,t)\\
&=C(n) a_1^n r_{\Rm}^n(p,t).
\end{aligned}
\end{equation}
Therefore, using the upper bound given by Corollary \ref{harnackcor}, we compute
\begin{equation}\label{eq.riemscaleupbound}
\begin{aligned}
\norm{r_{\Rm}^{-2}}_{\alpha,\beta,\mathcal{P}(p,t,a_1\,r_{\Rm})}&=\bigg( \int_{t-a_1^2 r_{\Rm}^2}^{t+a_1^2 r_{\Rm}^2 }{ \bigg( \int_{B_{g(t)}(p,a_1\,r_{\Rm}(p,t))}{\abs{r_{\Rm}^{-2}}^\alpha \, d\mu_s}\bigg)^{\beta/\alpha} ds} \bigg)^{1/\beta}\\
&\leq \bigg( \int_{t-a_1^2 r_{\Rm}^2}^{t+a_1^2 r_{\Rm}^2}{ \Big(4^{\alpha}C(n)a_1^n r_{\Rm}(p,t)^{n-2 \alpha} \Big)^{\beta/\alpha} ds} \bigg)^{1/\beta}\\
&=4\cdot 2^{1/\beta} C(n)^{1/\alpha} a_1^{n/\alpha+2/\beta} r_{\Rm}(p,t)^{n/\alpha + 2/\beta-2}\\
&= 4\cdot 2^{1/\beta} C(n)^{1/\alpha} a_1^2 \leq 8 C(n) a_1^2 =: C_3(n).
\end{aligned}
\end{equation}
Here we used in particular that $\frac{n}{\alpha}+\frac{2}{\beta}-2=0$ for an optimal pair. 

The proof of the opposite inequality follows a similar argument. Since $\abs{\Rm}$ is bounded by $r_{\Rm}^{-2}(p,t)$ on $\mathcal{P}(p,t,r_{\Rm})$ and $a_1\in(0,1)$, we can use the $\kappa$-noncollapsedness of the flow to obtain
\begin{equation*}
\mu_{g(t)}(B_{g(t)}(x,a_1\,r_{\Rm}(p,t)))\geq \kappa\, a_1^n\, r_{\Rm}^n(p,t).
\end{equation*}
Again, the evolution of the volume element under Ricci flow and the curvature bound in the cylinder considered yield the inequality
\begin{equation}\label{eq.volumelower}
\begin{aligned}
\mu_{g(s)}(B_{g(t)}(x,a_1\,r_{\Rm}(p,t))) &\geq e^{-n(n-1)r_{\Rm}^{-2}(p,t)(t-s)} \mu_{g(t)}(B_{g(t)}(x,a_1\,r_{\Rm}(p,t)))\\
&\geq \kappa\, e^{-n(n-1)a_1^2} a_1^n\, r_{\Rm}^n(p,t)\\
&=c(n,\kappa) a_1^n r_{\Rm}^n(p,t),
\end{aligned}
\end{equation}
for every $s \in (t-a_1^2 r^2_{\Rm}(p,t),t+a_1^2 r^2_{\Rm}(p,t))$. Corollary \ref{harnackcor} guarantees the lower bound on the integrand $r_{\Rm}^{-2}\geq \frac{1}{4} r_{\Rm}^{-2}(p,t)$ on the region $\mathcal{P}(p,t,a_1\,r_{\Rm})$. We can therefore follow \eqref{eq.riemscaleupbound}, reversing all inequalities except the very last one, to obtain
\begin{equation*}
\norm{r_{\Rm}^{-2}}_{\alpha,\beta,\mathcal{P}(p,t,a_1\,r_{\Rm})} \geq \frac{1}{4}\cdot 2^{1/\beta} c(n,\kappa)^{1/\alpha} a_1^2 \geq \frac{1}{4}c(n,\kappa)^{1/\alpha} a_1^2=: C_2(n,\kappa,\alpha).\qedhere
\end{equation*}
\end{proof}

The proof of Theorem \ref{thm.singintegralcharact} follows a very similar argument to what we have just seen. Notice, that since the integrand considered in the definition of the singular density function (Definition \ref{def.density}) is space-independent we simply get
\begin{equation}\label{eq.rewritenorm}
\Norm{\frac{1}{T-s}}_{\alpha,\beta,\mathcal{P}(p,t,a_1\,r_{\Rm})}= \bigg(\int_{t-a_1^2 r^2_{\Rm}}^{t+a_1^2 r^2_{\Rm}} \Big(\frac{1}{T-s}\Big)^\beta \mu_{g(s)}^{\beta/\alpha} \Big( B_{g(t)}(p,a_1\,r_{\Rm}(p,t)) \Big) ds \bigg)^{1/\beta}
\end{equation}
if $(t-a_1^2 r^2_{\Rm},t+a_1^2 r^2_{\Rm})\subseteq(0,T)$. We proceed with the proof.

\begin{proof}[Proof of Theorem \ref{thm.singintegralcharact}]
According to Perelman's noncollapsing theorem \cite{per}, there exists a constant $\kappa = \kappa(n,g(0),T)>0$ such that the Ricci flow in consideration is $\kappa$-noncollapsed at scale $\varrho=\sqrt{T}$. Hence, using \eqref{eq.volumeupper} and \eqref{eq.volumelower}, we obtain
\begin{equation}\label{eq.volumebounds}
c(n,\kappa)a_1^n r_{\Rm}^n(p,t) \leq \mu_{g(s)}\big( B_{g(t)}(p,a_1\,r_{\Rm}(p,t)) \big) \leq C(n)a_1^n r_{\Rm}^n(p,t),
\end{equation}
for some constants $c(n,\kappa)>0$ and $C(n)<\infty$, and where the lower bound requires $r_{\Rm}<\varrho$.

As in the proof of the previous theorem, if $p \in \Sigma$ we know from Corollary \ref{typeIoscicor} that $r^2_{\Rm}(p,t) \leq (T-t)$ and therefore \eqref{eq.volumebounds} holds for sufficiently large $t$ and we also have \eqref{eq.intinclusion} and \eqref{eq.nontruncated}. From \eqref{eq.intinclusion} we obtain in particular also
\begin{equation}\label{eq.Tsestimate}
\frac{2}{3(T-t)}\leq\frac{1}{T-s}\leq\frac{2}{T-t},\qquad \forall s \in (t-a_1^2 r^2_{\Rm}(p,t),t+a_1^2 r^2_{\Rm}(p,t)).
\end{equation}
Using \eqref{eq.rewritenorm}, we can therefore estimate
\begin{align*}
\Norm{\frac{1}{T-s}}_{\alpha,\beta,\mathcal{P}(p,t,a_1\,r_{\Rm})}&\leq \bigg(\int_{t-a_1^2 r^2_{\Rm}}^{t+a_1^2 r^2_{\Rm}} \Big(\frac{1}{T-s}\Big)^\beta \Big(C(n)a_1^n r_{\Rm}^n(p,t) \Big)^{\beta/\alpha} ds \bigg)^{\frac{1}{\beta}}\\
&\leq C(n)^{1/\alpha} (a_1\,r_{\Rm}(p,t))^{n/\alpha} \bigg(\int_{t-a_1^2 r^2_{\Rm}}^{t+a_1^2 r^2_{\Rm}} \Big(\frac{1}{T-s}\Big)^\beta ds \bigg)^{\frac{1}{\beta}}\\
&\leq C(n)^{1/\alpha} (a_1\,r_{\Rm}(p,t))^{n/\alpha}(2a^2_1 r^2_{\Rm}(p,t))^{1/\beta}\frac{2}{T-t}.
\end{align*}
Using that for an optimal pair $(\alpha,\beta)$ we have $\frac{n}{\alpha}+\frac{2}{\beta}=2$, we thus obtain
\begin{equation*}
\Norm{\frac{1}{T-s}}_{\alpha,\beta,\mathcal{P}(p,t,a_1\,r_{\Rm})}\leq 4C(n)^{1/\alpha}a^2_1\, \frac{r^2_{\Rm}(p,t)}{T-t}.
\end{equation*}
It is straightforward now to deduce from Theorem \ref{thm.singRiemannscale} that if $p \in \Sigma_I$, this is uniformly bounded from above by $4C(n)^{1/\alpha}a^2_1(\widetilde{c}_I)^{-1}$ for every $t$ close enough to $T$, whereas if $p \in \Sigma_{II}$ this gives $\Theta(p)=0$.  
Analogously, the lower bounds in \eqref{eq.volumebounds} and \eqref{eq.Tsestimate} yield
\begin{align*}
\Norm{\frac{1}{T-s}}_{\alpha,\beta,\mathcal{P}(p,t,a_1\,r_{\Rm})}&\geq \bigg(\int_{t-a_1^2 r^2_{\Rm}}^{t+a_1^2 r^2_{\Rm}} \Big(\frac{1}{T-s}\Big)^\beta \Big(c(n,\kappa)a_1^n r_{\Rm}^n(p,t) \Big)^{\beta/\alpha} ds \bigg)^{\frac{1}{\beta}}\\
&\geq c(n,\kappa)^{1/\alpha} (a_1\,r_{\Rm}(p,t))^{n/\alpha}(2a^2_1 r^2_{\Rm}(p,t))^{1/\beta}\frac{2}{3(T-t)}\\
&\geq\frac{2}{3}c(n,\kappa)^{1/\alpha}a^2_1\, \frac{r^2_{\Rm}(p,t)}{T-t}.
\end{align*}
If $p \in \Sigma_I$, then $\Theta(p)$ is bounded away from $0$ as $\liminf_{t\nearrow T} \frac{r^2_{\Rm}(p,t)}{T-t} \geq (\widetilde{C}_I)^{-1}$ by Theorem \ref{thm.singRiemannscale}. 

It remains to show that $\Theta(p)=\infty$ for regular points. For $p\in \mathfrak{Reg}$, Theorem \ref{thm.singRiemannscale} ensures the existence of a constant $C>0$ such that $r_{\Rm}^{-2}(p,t) \leq C^2$ for every $t \in [0,T]$. Thus for every $t \geq T- a_1^2/C^2$ we have $t+a^2_1r^2_{\Rm} > T$ and therefore
\begin{equation}\label{eq.Pinclusionregular}
\mathcal{P}(p,t,a_1\,r_{\Rm}(p,t)) \supseteq B_{g(t)}\Big(p,\frac{a_1}{C} \Big) \times (t,T).
\end{equation}
We can appeal to Perelman's noncollapsing theorem \cite{per} to obtain the existence of a constant $\kappa=\kappa(n,g(0),T,C)>0$ such that the Ricci flow in consideration is $\kappa$-noncollapsed at scale $\varrho=C^{-1}$. Therefore, since by (\ref{eq.Pinclusionregular})
\begin{equation*}
 \abs{\Rm} \leq r^{-2}_{\Rm}(p,t) \leq \frac{C^2}{a_1^2} \qquad \text{on} \qquad B_{g(t)}\Big(p,\frac{a_1}{C} \Big) \times (t,T),
\end{equation*}
we have for every $t \geq T-\frac{a_1^2}{C^2}$ and every $s \in (t,T)$ the uniform bound
\begin{equation*}
 \mu_{g(s)}\Big( B_{g(t)} \Big(p,\frac{a_1}{C} \Big) \Big) \geq c(n,\kappa) \frac{a_1^n}{C^n}>0.
\end{equation*}
Because $(T-s)^{-\beta} \notin L^1 \big(t,T\big)$ for every $t \geq T-\frac{a_1^2}{C^2}$, we obtain
\begin{align*}
\Norm{\frac{1}{T-s}}_{\alpha,\beta,\mathcal{P}(p,t,a_1\,r_{\Rm})}&\geq \Norm{\frac{1}{T-s}}_{\alpha,\beta,B_{g(t)}(p,\frac{a_1}{C}) \times (t,T)}\\
&\geq \Big(c(n,\kappa) \frac{a_1^n}{C^n}\Big)^{1/\alpha} \bigg(\int_{t}^{T} \frac{1}{(T-s)^{\beta}}ds \bigg)^{\frac{1}{\beta}}
=+\infty,
\end{align*}
from which we clearly see $\Theta(p)=+\infty$.
\end{proof}


\section{The Ricci Singular Sets}\label{sec.Ricci}

We begin this section with the definition of two different versions of Ricci scale, analogous to the definitions of the Riemann scales given in the introduction, Definition \ref{def.riemannscale}. Note that we always mark the fixed time-slice scales with a tilde in this article to distinguish them from the forwards-backwards scales that are mainly used in the local singularity analysis.

\begin{definition}[Ricci Scale]\label{def.ricciscale}
Let$(M,g(t))$ be a Ricci flow defined on $[0,T)$ and let $(p,t) \in M \times [0,T)$ be a space-time point. 
\begin{enumerate}[i)]
\item We define the \emph{Ricci scale} $r_{\Ric}(p,t)$ at (p,t) by
\begin{equation}\label{eq.Ricciscaledef}
r_{\Ric}(p,t) := \sup \{ r>0 \mid \abs{\Ric}<a_0(n)r^{-2} \text{ on } \mathcal{P}(p,t,r) \},
\end{equation}
where $a_0(n) = \sqrt{n}(n-1)$ as before, or equivalently by
\begin{equation}
r_{\Ric}(p,t) := \sup \{ r>0 \mid -(n-1)r^{-2}g < \Ric <(n-1)r^{-2}g \text{ on } \mathcal{P}(p,t,r) \}
\end{equation}
If $(M,g(t))$ is Ricci-flat for every $t\in[0,T)$, we set $r_{\Ric}(p,t)=+\infty$. Moreover, by slight abuse of notation, we may sometimes write $\mathcal{P}(p,t,r_{\Ric})$ for $\mathcal{P}(p,t,r_{\Ric}(p,t))$.
\item The \emph{time-slice Ricci scale} at $(p,t)$ is given by $\widetilde{r}_{\Ric}(p,t)=+\infty$ if the flow is Ricci flat, otherwise we set
\begin{equation}
\widetilde{r}_{\Ric}(p,t):= \sup \set{r>0 \mid \abs{\Ric}<a_0 r^{-2} \ \text{on} \ B_{g(t)}(p,r)},
\end{equation}
where $a_0=a_0(n):=\sqrt{n}(n-1)$.
\end{enumerate}
\end{definition}

\begin{remark}\label{rem.notourdefinition}
It might maybe seem more natural to define the Ricci scale as
\begin{equation}\label{eq.notourdefinition}
\sup \{ r>0 \mid \abs{\Ric}<r^{-2} \text{ on } \mathcal{P}(p,t,r) \}
\end{equation}
but our normalisation is more convenient for the purpose of this article for two main reasons. Firstly, because $\abs{\Rm}<r^{-2}$ implies 
\begin{equation}\label{eq.Riccibounds}
-(n-1)r^{-2}g < \Ric <(n-1)r^{-2}g,
\end{equation}
and thus $\abs{\Ric} <a_0r^{-2}$, we get the simple relation 
\begin{equation}\label{eq.RiccivsRmscale}
r_{\Ric}(p,t)\geq r_{\Rm}(p,t)
\end{equation}
for any space-time point $(p,t)$. Secondly, it is in fact exactly the property \eqref{eq.Riccibounds} which is used in a variety of proofs in the local singularity analysis and hence this normalisation allows a unified approach. (On the other hand, using \eqref{eq.notourdefinition} would allow to work with constants that do not depend on the dimension $n$, which could have advantages in other contexts.)
\end{remark}

\begin{remark}
From the Pseudolocality Proposition $3.2$ in \cite{bam} we see that $\widetilde{r}_{\Rm}$ and $r_{\Rm}$ are comparable for bounded scalar curvature Ricci flows. It is not clear whether a similar relation also holds for the Ricci scales, apart from the obvious estimate $\widetilde{r}_{\Ric}(p,t) \geq r_{\Ric}(p,t)$ that simply follows from the fact that the definition of $r_{\Ric}(p,t)$ requires a bound on a larger set.
\end{remark}

We first show that like the Riemann scale, also the Ricci scale is Lipschitz continuous in space and H\"older continuous in time, yielding in particular the following result.

\begin{theorem}[Lipschitz-H\"older Continuity of Ricci Scale]\label{thm.riccilipschitzholder}
Suppose $(M,g(t))$ is a complete $n$-dimensional Ricci flow defined on $[0,T)$. Then for any couple of space-time points $(p,t)$ and $(q,s)$ we have
\begin{equation}\label{eq.Riccilipschitzholder}
\abs{r_{\Ric}(p,t)-r_{\Ric}(q,s)}\leq \min \set{d_{g(t)}(p,q),d_{g(s)}(p,q)}+ C_1\abs{t-s}^{\frac{1}{2}},
\end{equation}
where $C_1=2\sqrt[4]{2/3} \sqrt{n-1}$ as in Theorem \ref{thm.lipschitzholder}.
\end{theorem}

\begin{proof}
We first note that if $(M,g(t))$ is a complete Ricci flow defined on $[0,T)$, then for any point $(p,t)\in M \times [0,T)$ we have the inclusion
\begin{equation}\label{eq.Ricciinclusion}
B_{g(t)}(p,r_{\Ric}(p,t)) \supseteq B_{g(s)}(p,r_{\Ric}(p,t)-C_1^2 r_{\Ric}^{-1}(p,t)\abs{s-t})
\end{equation}
for all $s\in(t-r_{\Ric}^2(p,t), t+r_{\Ric}^2(p,t))$. The proof of this inclusion is exactly the same as the one of Lemma \ref{shrinkingballslemma}, since it relies only on Ricci curvature bounds of the type \eqref{eq.Riccibounds}.

We can now deduce the Lipschitz-H\"older continuity exactly as in the proof of Theorem \ref{thm.lipschitzholder}. In fact, the only ingredient of the proof that does not rely on elementary estimates like the triangle inequality, is the use of \eqref{eq.Ricciinclusion} established above.
\end{proof}

A first corollary of this theorem gives upper and lower bounds for $r^{-2}_{\Ric}$ on a parabolic cylinder in terms of its value at the center of this cylinder. The proof is exactly the same as for the Riemann scale, exploiting the continuity from Theorem \ref{thm.riccilipschitzholder} above.
\begin{corollary}[Local Ricci Harnack-Type~Inequality]\label{ricciharnackcor}
Suppose $(M,g(t))$ is a complete Ricci flow on $[0,T)$, and let $(p,t)\in M \times (0,T)$. Then for $a_1=\frac{1}{2(1+C_1)}\in (0,1)$ as in Corollary \ref{harnackcor}, we have that
\begin{equation}\label{eq.Harnack}
\frac{1}{4} r^{-2}_{\Ric}(p,t)\leq r^{-2}_{\Ric}(\cdot,\cdot) \leq 4 r^{-2}_{\Ric}(p,t) \quad \text{on } \mathcal{P}(p,t,a_1r_{\Ric}(p,t)).
\end{equation}
\end{corollary}

Similar to Definition \ref{def.singpoints}, we can now define the following concepts of different types of \emph{Ricci singular sets}, depending only on the Ricci curvature rather than the full Riemannian curvature tensor.

\begin{definition}[Ricci Singular Points]\label{def.riccisingpoints}
Let $(M,g(t))$ be a Ricci flow on $[0,T)$, $T<\infty$ and assume that $(M,g(t))$ has bounded Ricci curvature on $M \times [0,t]$ for all $t\in[0,T)$. For any fixed $t\in[0,T)$, we consider again the parabolically rescaled Ricci flow $\widetilde{g}_t(s):= (T-t)^{-1}g(t+(T-t)s)$ defined for $s\in[-\frac{t}{T-t},1)$.
\begin{enumerate}[i)]
\item We say that a point $p \in M$ is a \emph{Ricci singular point} if for any neighbourhood $U$ of $p$, the Ricci curvature becomes unbounded on $U$ as $t$ approaches $T$. The \emph{Ricci singular set} $\Sigma^{\Ric}$ is the set of all such points and its complement is the \emph{Ricci regular set} $\mathfrak{Reg}^{\Ric}$.
\item We say that a point $p \in M$ is a \emph{Type~I Ricci singular point} if there exist constants $c_I,C_I, r_I>0$ such that we have
\begin{equation}\label{eq.riccitypeI}
a_0c_I < \limsup_{t\nearrow T} \sup_{B_{\widetilde{g}_t(0)}(p,r_I)\times(-r_I^2,r_I^2)} \,\abs{\Ric_{\widetilde{g}_t}}_{\widetilde{g}_t} \leq a_0C_I.
\end{equation}
We denote the set of such points by $\Sigma_I^{\Ric}$ and call it the \emph{Type~I Ricci singular set}.
\item We say that a point $p$ is a \emph{Type~II Ricci singular point} if for any $r>0$ we have
\begin{equation}\label{eq.riccitypeII}
\limsup_{t\nearrow T} \sup_{B_{\widetilde{g}_t(0)}(p,r)\times(-r^2,r^2)} \, \abs{\Ric_{\widetilde{g}_t}}_{\widetilde{g}_t}=\infty.
\end{equation}
We denote the set of such points by $\Sigma_{II}^{\Ric}$ and call it the \emph{Type~II Ricci singular set}.
\end{enumerate}
\end{definition}

We note that the upper bound in \eqref{eq.riccitypeI} is directly implied by the upper bound in \eqref{eq.typeI}. Obviously this is \emph{not} true for the respective lower bounds.

We immediately obtain the following alternative characterisations of the different Ricci singular sets which should be compared to their Riemann counterparts in Theorem \ref{thm.singRiemannscale}.

\begin{theorem}[Alternative Characterisation of Ricci Singular Sets]\label{thm.singRicciscale}
Let $(M,g(t))$ be a Ricci flow on $[0,T)$, $T<+\infty$ with bounded Ricci curvature on $M \times [0,t]$ for all $t \in [0,T)$. Let $\Sigma^{\Ric}$, $\Sigma_I^{\Ric}$, and $\Sigma_{II}^{\Ric}$ be given by Definition \ref{def.riccisingpoints}. Then
\begin{enumerate}[i)]
\item $p\in \Sigma^{\Ric}$ if and only if $\limsup_{t\nearrow T}\, r_{\Ric}^{-2}(p,t)=\infty$.
\item $p\in\Sigma_I^{\Ric}$ if and only if for some $0<\widetilde{c}_I, \widetilde{C}_I$ we have $\widetilde{c}_I < \limsup_{t\nearrow T}\, (T-t) r_{\Ric}^{-2}(p,t) \leq \widetilde{C}_I$.
\item $p\in\Sigma_{II}^{\Ric}$ if and only if $\limsup_{t\nearrow T}\, (T-t) r_{\Ric}^{-2}(p,t)=\infty$.
\end{enumerate}
\end{theorem}

\begin{proof}
The first part of the Theorem follows the same way as for the Riemann scale. The main point why this goes through is again that the proof relies on various applications of distance distortion estimates that only depend on the Ricci curvature bounds and not full Riemann bounds. The other parts of the proof can also be adopted almost verbatim, including the choices of all the constants -- for instance \eqref{eq.riccitypeI} implies
\begin{equation}\label{eq.riccialowerbound}
0 < \widetilde{c}_I:=\min \set{r_I^{-2},c_I} < (T-t_i)r^{-2}_{\Ric}(p,t_i) \leq \max \set{r_I^{-2}, C_I} =: \widetilde{C}_I.
\end{equation} 
This works because our normalisations in \eqref{eq.Ricciscaledef} and \eqref{eq.riccitypeI} agree.
\end{proof}

This then implies a non-oscillation result for the Ricci curvature along a Ricci flow, stating that it cannot oscillate between a Type~I rate and a lower rate arbitrarily close to the singular time. Once again, the proof is exactly the same as for the Riemann scale.
\begin{corollary}[Type~I non-Oscillation of Ricci scale]\label{RiccitypeIoscicor}
Suppose $(M,g(t))$ is a Ricci flow defined on a finite time interval $[0,T)$, with either bounded Ricci curvature on $M \times [0,t]$ for all $t\in[0,T)$ or complete time slices for all $t \in [0,T)$. Then $p \in \Sigma^{\Ric}$ if and only if
\begin{equation}
r^{-2}_{\Ric}(p,t) > \frac{1}{T-t}, \quad \forall t\in[0,T).
\end{equation}
\end{corollary}

\begin{remark}
This result may not yet be an improvement over Theorem $1$ in \cite{wan}, but it is its combination with Theorem \ref{thm.bridge} that gives a clear improvement from a global gap to a local one.
\end{remark}

As for the Riemann scale, combining Theorem \ref{thm.singRicciscale} with Corollary \ref{RiccitypeIoscicor}, we obtain the following decomposition of the Ricci singular set.

\begin{corollary}[Decomposition of the Ricci Singular Set]\label{cor.riccisingdecomposition}
Let $(M,g(t))$ be a Ricci flow on $[0,T)$, $T<+\infty$, with bounded Ricci curvature on $M\times [0,t]$ for every $t$ in $[0,T)$. Then $\Sigma^{\Ric}=\Sigma_I^{\Ric} \cup \Sigma_{II}^{\Ric}$.
\end{corollary}

As an immediate consequence of Corollary \ref{cor.riccisingdecomposition}, we also find the following result.
\begin{corollary}
Let $(M,g(t))$ be a Ricci flow on $[0,T)$, $T<+\infty$. Suppose the Ricci tensor satisfies a Type I bound. Then $\Sigma^{\Ric}=\Sigma_I^{\Ric}$.
\end{corollary}

We can now also prove the following $\e$-regularity type result.

\begin{theorem}[Integral Ricci Curvature Concentration]\label{ricciscaleintegralbounds}
Let $\kappa>0$ and let $(\alpha,\beta)$ be an optimal pair of integrability exponents in the sense of Definition \ref{def.optimal}. Then there exist constants $C_2=C_2(n,\kappa,\alpha)>0$ and $C_3=C_3(n)$ such that the following holds. Let $(M,g(t))$ be a complete Ricci flow defined on $[0,T)$, $T<\infty$, which is $\kappa$-non-local-collapsed on a scale $\varrho$ relative to the scalar curvature. Then for a space-time point $(p,t) \in \Sigma \times (T-\varrho^2,T)$, we have the integral bounds
\begin{equation}\label{eq.ricciintbounds}
C_2 \leq \norm{r_{\Ric}^{-2}}_{\alpha,\beta,\mathcal{P}(p,t,a_1\,r_{\Ric})} \leq C_3,
\end{equation}
where $a_1\in(0,1)$ is the constant from Corollary \ref{harnackcor}.
\end{theorem}

\begin{proof}
We note that, in addition to elementary estimates, the proof of Theorem \ref{riemannscaleintegralbounds} relies on the local Harnack type inequality for the Riemann scale (which holds with the same constants also for the Ricci scale), as well as distance and volume distortion estimates and the Bishop-Gromov inequality (all of which only rely on Ricci rather than full Riemann curvature bounds). Hence the proof can be adopted verbatim, simply changing every instance of $r_{\Rm}$ to $r_{\Ric}$. We leave it to the reader to check the details.
\end{proof}

In the remainder of this section, we focus on a localisation of the Sesum and Wang results, proving Theorem \ref{thm.bridge}. Our proof adapts some arguments from Theorem $2$ in Sesum \cite{ses} as well as ideas from Proposition $5.2$ in Hein-Naber \cite{hei}, combined with Theorem \ref{thm.singRicciscale} above.

\begin{proof}[Proof of Theorem \ref{thm.bridge}]
First of all, the inequality \eqref{eq.RiccivsRmscale} clearly implies, together with Theorem \ref{thm.singRicciscale}, that $\Sigma^{\Ric} \subseteq \Sigma$, so we only need to prove the opposite inclusion. In order to do so, we argue by contradiction and assume that there exists a point $p \in \Sigma$ such that $p \notin \Sigma^{\Ric}$. Using again Theorem \ref{thm.singRicciscale}, there exists a constant $\delta>0$ such that $r_{\Ric}(p,t)>\delta$ for every $t \in [0,T]$ while by Theorem \ref{thm.singRiemannscale} there exists a sequence of times $t_i \nearrow T$ so that $r_{\Rm}(p,t_i)\leq \frac{1}{i} \to 0$. For any $i \in \mathbb{N}$, let $q_i$ be a minimiser of the function $w_i$ defined by
\begin{equation}
w_i(q)=w_{(p,t_i)} (q) := \frac{r_{\Rm}(q,t_i)}{d_{g(t_i)}(q,\partial B_{g(t_i)}(p,\delta))},
\end{equation}
on the set $B_{g(t_i)}(p,\delta)$. We clearly have $w_i(q_i)\leq w_i(p)  = r_{\Rm}(p,t_i)/\delta \leq (i\delta)^{-1}$. As a consequence of the bounded curvature of the time-slices, we must also have that for every $i$, $r_{\Rm}(q_i,t_i)>0$, and therefore $w_i(q_i)>0$.

Set $r_i := r_{\Rm}(q_i,t_i)$ and consider the sequence of pointed rescaled Ricci flows $(M,g_i(t),q_i)$ defined by $g_i(t):=r_i^{-2}g(t_i+r_i^2 t)$ on $M\times [-r_i^{-2}t_i, r_i^{-2}(T-t_i))$. We first note that by definition of $q_i$, we have $r_i = r_{\Rm}(q_i,t_i) \leq r_{\Rm}(p,t_i)$ and thus, by Corollary \ref{typeIoscicor}, $r_i^{-2}(T-t_i) \geq 1$, hence the flows $g_i(t)$ exist at least for times $t\in[-1, 1)$. By definition, they satisfy $r_{\Rm_{i}}(q_i,0)=1$ for every $i$ and by the scaling properties of the distance
\begin{equation*}
d_i := \frac{1}{2} d_{g_i(0)}(q_i, \partial B_{g_i(0)}(p,\delta r_i^{-1})) = \frac{1}{2 w_{(p,0)}(q_i)} \ge \frac{\delta i}{2} \to + \infty.
\end{equation*}
Notice that by definition of $d_i$ we obtain $B_{g_i(0)}(q_i,d_i)\subseteq B_{g_i(0)}(p,\delta r_i^{-1})= B_{g(t_i)}(p,\delta)$. Since for every $q \in B_{g_i(0)}(q_i,d_i)$ its $g_i(0)$-distance to $\partial B_{g_i(0)}(p,\delta r_i^{-1})$ is at least $d_i$, we deduce from the minimising property of $q_i$ that
\begin{equation*}
\frac{1}{2 d_i} =w_{(p,t_i)}(q_i) \leq w_{(p,t_i)}(q) \leq \frac{r_{\Rm_i}(q,0)}{d_i} \iff r_{\Rm_i}(q,0) \ge \frac{1}{2}.
\end{equation*}
Perelman's non local-collapsing theorem \cite{per} applied to any of the cylinders $\mathcal{P}(q,0,r_{\Rm_i}(q,0))$ with $q \in B_{g_i(0)}(q_i,d_i)$ as before, guarantees the existence of a uniform injectivity radius lower bound. Therefore, we can apply Topping's compactness theorem \cite[Theorem 1.6]{top} to extract a pointed smooth Cheeger-Gromov limit Ricci flow $(M_\infty,g_\infty(t),q_\infty)$, defined and complete in $M_\infty \times (-\frac{1}{4},\frac{1}{4})$. To do so, we just need to check that for every $r>0$ there exists (in Topping's notation) some $K(r)\in \mathbb{N}$ such that for every $i \ge K(r)$ the curvature is uniformly bounded by $M$ on $B_{g_i(0)}(q_i,r)$. But this is obviously true with $M=4$ and $K(r)$ such that $d_i \ge r$ for $i \ge K(r)$.

This limit flow inherits several properties. First of all, $r_{\Rm_\infty}(q_\infty,0)=1$ and $\abs{\Rm_\infty}\leq 4$ on $M_\infty \times (-\frac{1}{4},\frac{1}{4})$. Secondly, we deduce from the inclusion $B_{g_i(0)}(q_i,d_i) \subseteq B_{g(t_i)}(p,\delta)$ that $\abs{\Ric_i} \le a_0r_i^2 \delta^{-2}$, so the limit flow must satisfy $\Ric_\infty \equiv 0$, i.e. $g_{\infty}(t) \equiv g_{\infty}$ is a static Ricci flat metric. We have therefore reconducted the study to a situation similar to the case treated by Sesum \cite{ses}.

For any $r>0$, the smooth Cheeger-Gromov convergence ensures
\begin{equation*}
\frac{\mu_{g_\infty}(B_{g_\infty}(q_\infty,r))}{\omega_n r^n}=\lim_{i \rightarrow +\infty} \frac{\mu_{g_i(0)}(B_{g_i(0)}(q_i,r))}{\omega_n r^n}=\lim_{i \rightarrow +\infty} \frac{\mu_{g(t_i)}(B_{g(t_i)}(q_i,r r_i))}{\omega_n (r r_i)^n}.
\end{equation*}
Given $\delta>0$ as above, there exists $i_0(\delta,r)$ such that $r r_i < \delta$ and $t_i+\delta^2 >T$ for every $i \ge i_0$, therefore we can appeal to the bound $\abs{\Ric} \le a_0\delta^{-2}$ to use the multiplicative distance distortion estimates in $B_{g(t_i)}(p,r r_i) \times [t_{i_0}-\delta^2,T)$ for every $i \geq i_0$. For any $\e>0$ and a possibly even larger $i_1(\delta,\e)$, we have the following two conditions satisfied for all $i\geq i_1$ 
\begin{equation*}
(t_i-t_{i_1})\le a_0^{-1}\delta^2 \e,
\end{equation*}
and
\begin{equation*}
\frac{\mu_{g_\infty}(B_{g_\infty}(q_\infty,r))}{\omega_n r^n}\ge \frac{\mu_{g(t_i)}(B_{g(t_i)}(q_i,r r_i))}{\omega_n (r r_i)^n} -\frac{\e}{2}.
\end{equation*}
Using the distortion estimate we see that
\begin{equation*}
 B_{g(t_i)}(q_i,r r_i) \supseteq B_{g(t_{i_1})}(q_i,r r_i e^{-a_0\delta^{-2}(t_i-t_{i_1})})\supseteq B_{g(t_{i_1})}(q_i,r r_i c_{\e}).
\end{equation*}
Here $c_{\e}=e^{-\e} \rightarrow 1$ as $\e \rightarrow 0$.
Therefore, by the bound on the scalar curvature and the evolution equation for the volume element we obtain
\begin{align*}
\mu_{g(t_i)}(B_{g(t_{i_1})}(q_i,r r_i c_{\e})) &\ge e^{-\sqrt{n}\,a_0 \delta^{-2}(t_i-t_{i_1})} \mu_{g(t_{i_1})}(B_{g(t_{i_1})}(q_i,r r_i c_{\e})) \\
&\ge c_{\e}^{\sqrt{n}} \mu_{g(t_{i_1})}(B_{g(t_{i_1})}(q_i,r r_i c_{\e})).
\end{align*}
We have arrived at
\begin{equation*}
\frac{\mu_{g_\infty}(B_{g_\infty}(q_\infty,r))}{\omega_n r^n}\ge c_{\e}^{\sqrt{n}}\, \frac{\mu_{g(t_{i_1})}(B_{g(t_{i_1})}(q_i,c_{\e}r r_i))}{\omega_n (r r_i)^n} -\frac{\e}{2}.
\end{equation*}
The key feature of this inequality is that the right hand side has a controllable dependence on $i$. Arguing as in \cite{ses}, using the fact that we have full control on the geometry of $g(t_{i_1})$ by the bounded curvature assumption and Shi's estimates as well as the fact that $r_i\to 0$, we obtain
\begin{equation*}
 \mu_{g(t_{i_1})}(B_{g(t_{i_1})}(q_i,c_{\e}r r_i)) \ge  \omega_n (c_{\e}r r_i)^n \Big(1-\frac{\e}{2}\Big)
\end{equation*}
for sufficiently large $i \ge i_2(\e,g_{i_1},n)$. In particular, sending $i$ to infinity, we obtain
\begin{equation}
\frac{\mu_{g_\infty}(B_{g_\infty}(q_\infty,r))}{\omega_n r^n}\ge \Big(1-\frac{\e}{2}\Big) c_{\e}^{n+\sqrt{n}} -\frac{\e}{2}.
\end{equation}
We can now let $\e$ go to $0$ to find
\begin{equation}
\frac{\mu_{g_\infty}(B_{g_\infty}(q_\infty,r))}{\omega_n r^n}\ge 1.
\end{equation}
On the other hand, since $g_\infty$ is Ricci flat, by Bishop-Gromov inequality we obtain
\begin{equation}
\frac{\mu_{g_\infty}(B_{g_\infty}(q_\infty,r))}{\omega_n r^n}\le 1,
\end{equation}
and hence we are in the equality case where $B_{g_\infty}(q_\infty,r)$ has exactly Euclidean volume growth. By the rigidity statement in Bishop-Gromov's inequality, we deduce that $g_\infty$ is flat, in contradiction with $r_{\Rm_\infty}(q_\infty)=1$. This means that $\Sigma \subseteq \Sigma^{\Ric}$ as claimed.
\end{proof}

It is now easy to prove Corollary \ref{cor.bridge}.
\begin{proof}
Fix a point $p \in \Sigma_{II}^{\Ric}$. Then by Theorem \ref{thm.singRicciscale} we know that
\begin{equation}
\limsup_{t \nearrow T} (T-t) r_{\Ric}^{-2}(p,t)=+ \infty,
\end{equation}
and hence by the inequality $r_{\Rm}^{-2}(p,t) \ge r_{\Ric}^{-2}(p,t)$ also
\begin{equation}
\limsup_{t \nearrow T} (T-t) r_{\Rm}^{-2}(p,t)=+ \infty.
\end{equation}
Theorem \ref{thm.singRiemannscale} $iii)$ guarantees that $p \in \Sigma_{II}$. The other inclusion is easily deduced from this one, Theorem \ref{thm.singdecomposition}, Corollary \ref{cor.riccisingdecomposition}, and Theorem \ref{thm.bridge}.
\end{proof}

Arguing as we did in Proposition \ref{prop.singRiemannscale} the reader can easily verify the following result.
\begin{proposition}[Characterisation of Ricci Singular Set using Fixed Time Slice Scale]\label{prop.singRicciscale}
Let $(M,g(t))$ be a complete Ricci flow on $[0,T)$, $T<\infty$ with bounded Ricci curvature on $M\times [0,t]$ for every $t \in [0,T)$. Then $p \in \Sigma^{\Ric}$ if and only if $\liminf_{t \nearrow T}\widetilde{r}_{\Ric}(p,t)=0$.
\end{proposition}
Together with our Theorem \ref{thm.bridge} and Proposition \ref{prop.singRiemannscale} , we deduce the following corollary.
\begin{corollary}
Let $(M,g(t))$ be a Ricci flow defined on $[0,T)$, $T<\infty$, such that $(M,g(t))$ is complete and has bounded curvature on $M\times [0,t]$ for every $t \in [0,T)$. Suppose the initial time slice satisfies $\inj(M,g(0))>0$. Then for any point $p \in M$, $\liminf_{t \nearrow T}\widetilde{r}_{\Rm}(p,t)=0$ if and only if $\liminf_{t \nearrow T}\widetilde{r}_{\Ric}(p,t)=0$.
\end{corollary}

\section{Applications to Bounded Scalar Curvature Ricci Flows}\label{sec.scalar2}

The aim of this section is to discuss applications of our local theory to Ricci flows with bounded scalar curvature. We first recall several results from Bamler \cite{bam} and Bamler-Zhang \cite{bamz}. 

We first point out that in our discussion about noncollapsing relative to the scalar curvature (see Definition \ref{def.noncollapsed} and the paragraph below), we worked with an initial injectivity radius bound and an implicit initial Ricci curvature lower bound. Alternatively, we could have assumed a lower bound on Perelman's entropy of the initial metric, $\nu_0:=\nu[g(0),2T]>-A$ to obtain these noncollapsing results. This is the condition present in Bamler's work \cite{bam}, but we will rewrite the results we need from his paper using injectivity radius bounds instead; we start by recalling his extremely powerful estimate of the volume of sublevel sets of $\widetilde{r}_{\Rm}$ for Ricci flows with bounded scalar curvature.

\begin{theorem}[Proposition $6.4$ of \cite{bam}]\label{th.proposition64bamler}
For any $n \in \mathbb{N}$, $T<+\infty$, $R_0>0$, $i_0>0$, $k_0>0$ and $d' \in (0,4)$ there exist constant $E'=E'(n,R_0,i_0,k_0,d')$ and $t_0=t_0(n,R_0,T)<T$ such that the following statement holds. Let $(M^n,g(t))$ be a closed Ricci flow defined on $[0,T)$ and satisfying $\inj(M,g(0))> i_0$ and $\Ric_{g(0)} \ge -(n-1) k_0 g(0)$. Assume that the scalar curvature is uniformly bounded, $\abs{\Sc} \le n (n-1)R_0$. Then for any $p \in M$, $t \in [t_0,T)$ and $r,s \in (0,1)$, we have
\begin{equation}\label{eq.64bamler}
\mu_{g(t)}(\set{p \mid \widetilde{r}_{\Rm}(p,t)<sr} \cap B_{g(t)}(p,r))< E' s^{d'} r^n. 
\end{equation}
\end{theorem}
Bamler uses this estimate to prove a codimension four estimate of the singular set (after passing to a weak limit). We will use it for our codimension eight result, Theorem \ref{th.codimension8}. Let us remark that in its original wording, Proposition $6.4$ of \cite{bam} requires $\abs{\Sc}\leq 1$ on a time interval $[-2,0]$ (yielding $t_0=-1$), but the result as stated above can be obtained from this by parabolic rescaling. The proof of the above theorem is extremely involved and occupies most of \cite{bam}. It relies on a detailed analysis on the geometry of Perelman's reduced length, which allows the author to improve further and further the bound on the volume of high curvature regions. Luckily, we do not need to modify any of these estimates but manage to use Theorem \ref{th.proposition64bamler} as a ``black box''.

The theorem has the following corollary.

\begin{corollary}[Theorem $1.7$ of \cite{bam}]\label{cor.theorem17bamler}
For any $n \in \mathbb{N}$, $T<+\infty$, $R_0>0$, $i_0>0$, $k_0>0$ and $\e>0$, there exist constants $F=F(n,R_0,i_0,k_0,\e)$ and $t_0=t_0(n,R_0,T)<T$ as in the theorem above, such that the following statement holds. Let $(M^n,g(t))$ be a closed Ricci flow defined on $[0,T)$, $T<+\infty$ and satisfying $\inj(M,g(0))> i_0$ and $\Ric_{g(0)} \ge -(n-1) k_0 g(0)$. Assume that the scalar curvature is uniformly bounded, $\abs{\Sc} \le n (n-1)R_0$. Then for every $p\in M$, $r' \in (0,1)$ and $t \in [t_0,T)$, we have
\begin{equation}
\int_{B_{g(t)}(p,r')} \widetilde{r}_{\Rm}^{\,-\alpha}(\cdot,t) \, d\mu_{g(t)} \le F (r')^{n-4+2\e},
\end{equation}
where $\alpha:=4-2\e$.
\end{corollary}

As this result follows directly from Theorem \ref{th.proposition64bamler}, we give a quick sketch here in an attempt to make this article more self-contained.

\begin{proof}
For $\alpha=4-2\e$, set $d':=4-\e$ and pick $E'$ and $t_0$ as in Theorem \ref{th.proposition64bamler}. Similar to the noncollapsing volume bound, the assumptions on our flow ensure a noninflation result, meaning that for some constant $K_1=K_1(n,i_0,k_0,R_0)$ we have 
\begin{equation}
\mu_{g(t)} (B_{g(t)}(p,r)) \leq K_1 r^n
\end{equation} 
for all $t \in [t_0,T)$, see e.g. \cite{z12}. The dependency of the constant $K_1$ on the initial metric $g(0)$ is not explicitly related to $i_0$ and $k_0$ in \cite{z12}, but it is shown that it depends only on the initial $\mathcal{F}$-entropy and initial Sobolev constants bounds, which can be both bounded in terms of $i_0$ and $k_0$ as remarked by the same author in \cite{z07}, see the remark after Theorem $1.1$ there.

We can therefore estimate
\begin{align*}
I(r') &:= \int_{B_{g(t)}(p,r')} \widetilde{r}_{\Rm}^{\,-\alpha}(q,t) \, d\mu_{g(t)}(q)\\
&\leq \int_{B_{g(t)}(p,r')} \bigg( \int_0^{(r')^{-\alpha}} 1\, dy + \int_{(r')^{-\alpha}}^{\infty} \chi_{\{y<\widetilde{r}_{\Rm}^{\,-\alpha}(q,t)\}}(q,y)\, dy\bigg) d\mu_{g(t)}(q)\\
&= (r')^{-\alpha} \mu_{g(t)} (B_{g(t)}(p,r')) + \int_{(r')^{-\alpha}}^{\infty} \mu_{g(t)}(\{q \mid \widetilde{r}_{\Rm}(q,t) < y^{-1/\alpha}\} \cap B_{g(t)}(p,r'))\, dy\\
&\leq K_1 (r')^{n-4+2\e} + \int_{(r')^{-\alpha}}^{\infty} E' y^{-d'/\alpha} (r')^{-d'} (r')^n \, dy\\
&\leq K_1 (r')^{n-4+2\e} + E' (r')^{n-4+\e} \int_{(r')^{-\alpha}}^{\infty} y^{-d'/\alpha} dy\\
&= F (r')^{n-4+2\e},
\end{align*}
where $F:=K_1 + E' \e/\alpha$. Here, we have used \eqref{eq.64bamler} with $s:=y^{-1/\alpha}/r'$ and the fact that the very last integral is equal to $(r')^{\e}\cdot \e/\alpha$ for the choices of $d'$ and $\alpha$ as above.
\end{proof}

As a last ingredient for our results, we recall the following lemma, essentially due to Wang \cite{wan}, though our version resembles more the one from Bamler and Zhang \cite{bamz}, clarifying in which terms the square root of the Riemann curvature controls the Ricci curvature on a flow with bounded scalar curvature.

\begin{lemma}[Lemma $6.1$ in \cite{bamz}]\label{th.lemma6.1BZ}
For any $n \in \mathbb{N}$, $R_0>0$, $i_0>0$ and $k_0>0$, there exists a constant $C_4=C_4(n,R_0,i_0,k_0)$, such that the following statement holds. Let $(M^n,g(t))$ be a closed Ricci flow defined on $[0,T)$, $T<+\infty$ and satisfying $\inj(M,g(0))> i_0$ and $\Ric_{g(0)} \ge -(n-1) k_0 g(0)$. Suppose that for some $(p_0,t_0) \in M \times (0,T)$ and $r_0 \in (0,\min \{ R_0^{-1/2}, \sqrt{t_0}, 1 \})$ we have $B_{g(t_0)}(p_0,r_0) \subset \subset M$. If $\abs{\Sc} \le n(n-1)R_0$ and $\abs{\Rm} \le r_0^{-2}$ on the backward cylinder $\mathcal{P}^{-}(p_0,t_0,r_0):=B_{g(t_0)}(p_0,r_0)\times (t_0-r_0^2,t_0]$, then we have the bound $\abs{\Ric(p_0,t_0)} \le C_4 a_0 r_0^{-1}$.
\end{lemma}

As this result is essential for our applications of the local theory to bounded scalar curvature Ricci flows, we recall the proof here, closely following the one given in \cite{bamz}.

\begin{proof}
By parabolically rescaling the flow with factor $r_0^{-2}$ we may assume $\abs{\Rm} \le 1$ on the cylinder $\mathcal{P}^{-}(p_0,t_0,1)$ and $\abs{\Sc} \le A := n(n-1)R_0 r_0^2$. We need to show the existence of a constant $D(n,R_0,i_0,k_0)$ such that $\abs{\Ric(p_0,t_0)} \le D \sqrt{A}$, yielding the claim for $C_4= D\sqrt{R_0}/\sqrt{n-1}$.

Using Shi's estimates we see that for universal constants $\{C_m\}_{m \in \mathbb{N}}$ we have
\begin{equation}\label{eq.smoothbounds}
|\nabla^m \Rm| \le C_m \quad \text{on } \mathcal{P}^{-}\big(x_0,t_0,\tfrac{1}{2}\big).
\end{equation}
Notice that by our assumption, we have that $\mathcal{P}^{-}(x_0,t_0,a_1 r_0) \subset \subset M \times (0,T)$, where $a_1$ is the constant defined in Section \ref{sec.pointwise}. Arguing as in \cite{bamz} we deduce the existence of a universal constant $b_1 \le \min \{ a_1, 1/2\}$ such that the map $\exp_{p_0}\colon B(0,b_1)\subset \mathbb{R}^n \to M$ given by the $g(t_0)-$exponential map centered at $p_0$ is injective, and the pull-back Ricci flow $\widetilde{g}(t) := \exp_{p_0}^*g(t)$ defined on $B(0,b_1) \times [t_0-b_1^2,t_0]$ inherits smooth bounds by (\ref{eq.smoothbounds}), with the metrics $\widetilde{g}(t)$ being $2$-Lipschitz equivalent to the Euclidean metric for every $t \in [t_0-b_1^2,t_0]$. Note that the bounds on the curvatures are preserved by the pull-back (up to universal constants). From now on we will focus on the metrics $\widetilde{g}(t)$, and drop the tilde to simplify the notation. Fix a cut-off function $\phi \in C^{\infty}_0(B(0,b_1))$ such that $\phi \in [0,1]$ and $\phi\equiv 1$ on $B(0,\frac{b_1}{2})$. This can be done in such a way that $\abs{\partial \phi}, \abs{\partial^2 \phi} \le E_1$ for some universal constant $E_1$. Therefore, for some other universal $E_2$ we have $|\Delta_{g(t)} \phi|\le E_2$ on $B(0,b_1) \times [t_0-b_1^2,t_0]$.
Testing the evolution equation for the scalar curvature against $\phi$, and integrating by parts we obtain
\begin{align*}
\bigg| \partial_t \int_{B(0,b_1)} \Sc(\cdot,t) \phi\, d\mu_{g(t)} &-\int_{B(0,b_1)} 2|\Ric(\cdot,t)|^2 \phi\, d\mu_{g(t)}\bigg|\\
&=\bigg| -\int_{B(0,b_1)} \Sc(\cdot,t)^2 \phi\, d\mu_{g(t)}+\int_{B(0,b_1)} \Delta_{g(t)} \Sc(\cdot,t) \phi\, d\mu_{g(t)}\bigg|\\
&=\bigg| -\int_{B(0,b_1)} \Sc(\cdot,t)^2 \phi\, d\mu_{g(t)}+\int_{B(0,b_1)} \Delta_{g(t)} \phi \Sc(\cdot,t) d\mu_{g(t)}\bigg|.
\end{align*}
The inequality $\abs{\Sc}\le A \le n(n-1)R_0$ implies for every $t \in [t_0-b_1^2,t_0]$ the non-inflating property $\mu_{g(t)}(B(0,b_1))\le K_1 b_1^n := E_3$ for some $K_1=K_1(n,i_0,k_0,R_0)$ by \cite{z12}. Moreover, by the discussion above, the Laplace term is controlled, so we get the upper bound for the right-hand side
\begin{align*}
\bigg| -\int_{B(0,b_1)} \Sc(\cdot,t)^2 \phi\, d\mu_{g(t)} &+\int_{B(0,b_1)} \Delta_{g(t)} \phi \Sc(\cdot,t) d\mu_{g(t)}\bigg|\\
&\le E_3 A^2+E_3 E_2 A \le (n(n-1)R_0E_3+E_2 E_3) A=:E_4 A.
\end{align*}
The constant $E_4$ depends only on $n$, $i_0$, $k_0$ and $R_0$. Integrating this inequality in time we deduce
\begin{align*}
\norm{\Ric}_{2,\mathcal{P}^{-}(0,b_1/2)}^2 &\le \int_{t_0-b_1^2}^{t_0} \int_{B(0,b_1)} \abs{\Ric(\cdot,t)}^2 \phi\, d\mu_{g(t)} dt \\
&\le \int_{B(0,b_1)} \abs{\Sc(\cdot,t_0)} \phi\, d\mu_{g(t_0)}+\int_{B(0,b_1)} \abs{\Sc(\cdot,t_0-b_1^2)} \phi\, d\mu_{g(t_0-b_1^2)} +b_1^2 E_4 A\\
&\le 2E_3 A+b_1^2 E_4A=:E_5 A.
\end{align*}
Recall that the Ricci tensor solves the parabolic system
\begin{equation}
(\partial_t+\Delta_{g(t)}-2\Rm)\Ric=0,
\end{equation}
which can be interpreted as linear in $\Ric$, with coefficients universally bounded in every $C^m$ norm. Thus, the standard parabolic theory ensures the existence of a universal constant $E_6$ such that, once we set $b_2:=b_1/4$, we have
\begin{equation}
\abs{\Ric(0,t_0)} \leq \norm{\Ric}_{\infty,\mathcal{P}^{-}(0,b_2)} \le E_6 \norm{\Ric}_{2,\mathcal{P}^{-}(0,b_1/2)} \le E_6 \sqrt{E_5 A} =: D \sqrt{A},
\end{equation}
concluding the proof.
\end{proof}

The proof above shows that one can extend the Ricci bounds to the backward cylinder $\mathcal{P}^{-}(p_0,t_0,b_2r_0)$ for a universal constant $b_2\in(0,1)$. A backward-forward analogue of this result would be more in line with the results in the rest of this paper and it is indeed also possible. Nevertheless, in each case, the bound obtained is not strong enough to get the desired estimate $\widetilde{r}_{\Ric}^{\,2} \gtrsim \widetilde{r}_{\Rm}$, because this would require a Ricci bound on the (bigger) scale $\sqrt{r_0}$. This issue comes directly from the proof as given above, since one can appeal to linear parabolic regularity theory only on the scale $r_0$. A proof of a global relation $\widetilde{r}_{\Ric}^{\, 2} \gtrsim \widetilde{r}_{\Rm}$ would have extremely interesting consequences, for example it would imply $\Sigma_I=\emptyset$, but without extra assumptions, it currently seems out of reach.
 
Here, we overcome this difficulty by considering only points where the Ricci curvature is almost maximal or well behaved points as defined in the introduction, so that any bound on the Ricci curvature that we find extends naturally to a bound on the Ricci scale.

\begin{theorem}[Quadratic Scale Comparison at Certain Points]\label{thm.squareroot}
For any $n \in \mathbb{N}$, $T<+\infty$, $R_0>0$, $i_0>0$ and $k_0>0$ there exists a constant $a_3=a_3(n,R_0,i_0,k_0,T)$ such that the following holds.
Let $(M^n,g(t))$ be a closed Ricci flow defined on a finite time interval $[0,T)$, satisfying $\inj(M,g(0))>i_0$ and $\Ric_{g(0)} \ge -(n-1) k_0 g(0)$. Assume that $\abs{\Sc}\le n(n-1)R_0$ on $M \times [0,T)$. Then we deduce the two properties below.
\begin{enumerate}[i)]
\item Let $(p_0,t_0) \in M\times(0,T)$ be such that $\frac{1}{2}\sup_M \abs{\Ric}(\cdot,t_0) \leq \abs{\Ric}(p_0,t_0) =: a_0 r_0^{-2}$. If $r_0 \le \min \{R_0^{-1/2},\sqrt{t_0},1\}$, we have $\widetilde{r}^{\, 2}_{\Ric}(p_0,t_0) \ge a_3 \widetilde{r}_{\Rm}(p_0,t_0)$.
\item Let $(p_0,t_0) \in G_{\delta,t_1}$ for some $\delta \in (0,1)$, where $G_\delta=G_{\delta,t_1}$ denotes the set of well behaved points for some $t_1\leq t_0$, see \eqref{eq.wellbehaved}. Then if $r_{\Ric}(p_0,t_0) \le \min \{R_0^{-1/2},\sqrt{t_0},1\}$, we have $r^{2}_{\Ric}(p_0,t_0) \ge \delta a_3 r_{\Rm}(p_0,t_0)$.
\end{enumerate}
\end{theorem}

\begin{proof}
The proof is essentially the same for both cases.
\begin{enumerate}[i)]
\item Fix any such space-time point $(p_0,t_0)$ with almost maximal Ricci curvature. From the definition of the time-slice Ricci scale, it is clear that $a_0\widetilde{r}_{\Ric}^{\,-2}(p_0,t_0) \geq \abs{\Ric}(p_0,t_0) = a_0r_0^{-2}$. Moreover, since $p_0$ has almost maximal Ricci curvature, $\abs{\Ric}(\cdot,t_0) \leq 2\abs{\Ric}(p_0,t_0)$ everywhere and thus $a_0\widetilde{r}_{\Ric}^{\,-2}(p_0,t_0) \leq 2\abs{\Ric}(p_0,t_0) = 2a_0r_0^{-2}$. Together, we have 
\begin{equation}\label{eq.ricciscalesandwich}
r_0^{-2} \leq \widetilde{r}_{\Ric}^{\,-2}(p_0,t_0) \leq 2r_0^{-2}.
\end{equation}
By definition of the time-slice Riemann scale, we have $\abs{\Rm} \le \widetilde{r}_{\Rm}^{\, -2}(p_0,t_0)$ on the ball $B_{g(t_0)}(p_0,\widetilde{r}_{\Rm}(p_0,t_0))$. The Backward Pseudolocality Theorem of Bamler-Zhang (Theorem $1.5$ in \cite{bamz}) applied to the scale $r_1:=\widetilde{r}_{\Rm}(p_0,t_0) \le \widetilde{r}_{\Ric}(p_0,t_0) \leq r_0 \le \sqrt{t_0}$ yields the existence of constants $\e$ and $K$ (without loss of generality $K \le \e^{-2}$) depending only on $n$, $T$, $i_0$ and $k_0$ such that we have $\abs{\Rm} \le K r_1^{-2}\le (\e r_1)^{-2}$ on the backwards parabolic cylinder $\mathcal{P}^{-}(p_0,t_0,\e r_1)$. Using Lemma \ref{th.lemma6.1BZ}, we get a constant $C_4=C_4(n,R_0,i_0,k_0)$ such that $a_0 r_0^{-2}=\abs{\Ric}(p_0,t_0) \le C_4a_0(\e r_1)^{-1}$ or equivalently
\begin{equation*}
\widetilde{r}_{\Ric}^{\,-2}(p_0,t_0) \leq 2r_0^{-2} \leq 2C_4(\e r_1)^{-1} = \frac{2C_4}{\e}\, \widetilde{r}_{\Rm}^{\, -1}(p_0,t_0).
\end{equation*} 
The conclusion now follows for $a_3\leq \frac{\e}{2C_4}$.
\item Fix a well behaved point $(p_0,t_0) \in G_\delta$. Denoting again $\abs{\Ric}(p_0,t_0) =: a_0 r_0^{-2}$, we obtain the following analogue of equation (\ref{eq.ricciscalesandwich}) 
\begin{equation}
r_0^{-2} \leq r_{\Ric}^{-2}(p_0,t_0) \leq \delta^{-1}r_0^{-2}.
\end{equation}
Since $r_{\Rm}(p_0,t_0)\leq r_{\Ric}(p_0,t_0) \le \min \{R_0^{-1/2},\sqrt{t_0},1\}$, and we have, in particular, $\abs{\Rm} \le r_{\Rm}^{-2}(p_0,t_0)$ on the (backwards) parabolic cylinder $\mathcal{P}^{-}(p_0,t_0,r_{\Rm})$, we can apply again Lemma \ref{th.lemma6.1BZ} to conclude $a_0 r_0^{-2}=\abs{\Ric}(p_0,t_0) \le C_4a_0(r_{\Rm}(p_0,t_0))^{-1}$. This means that
\begin{equation*}
r_{\Ric}^{-2}(p_0,t_0) \leq \delta^{-1} r_0^{-2} \leq \delta^{-1}C_4 r_{\Rm}^{-1}(p_0,t_0)),
\end{equation*} 
and hence the conclusion for $a_3\leq \frac{1}{C_4}$.\qedhere
\end{enumerate}
\end{proof}

Next, we prove the following integral Ricci curvature concentration result, which is the time-slice analogue of Theorem \ref{ricciscaleintegralbounds}.

\begin{lemma}[Time-Slice Ricci Scale Integral Concentration]\label{lemma.L8lowerbound}
Let $(M^n,g(t))$ be a closed $n$-dimensional Ricci flow defined on a finite time interval $[0,T)$ with $\inj(M,g(0))>i_0$ and $\Ric_{g(0)}\ge -(n-1)k_0 g(0)$ for some constants $i_0,k_0>0$. Assume there exists $R_0<\infty$ such that $\abs{\Sc}\le n(n-1)R_0$ on $M \times [0,T)$, and let $\alpha \ge \frac{n}{2}$. Then for any $\delta\in (0,\frac{1}{2})$, there exists $C_5=C_5(n,R_0,T,i_0,k_0,\alpha,\delta)>0$ such that for any $(p_0,t_0) \in M \times [0,T)$ with $\widetilde{r}_{\Ric}(p_0,t_0)\le \min \set{R_0^{-1/2},\sqrt{T}}$, we can bound
\begin{equation}
\int_{B_{g(t_0)}(p_0,\delta\widetilde{r}_{\Ric}(p_0,t_0))} \widetilde{r}_{\Ric}^{\,-2 \alpha} \, d\mu_{g(t_0)} \ge C_5>0.
\end{equation}
\end{lemma}

\begin{proof}
We first note that while the fixed time slice scales might not satisfy nice continuity properties in time, their spatial continuity is well understood. In fact, for any $t \in [0,T)$, the functions $\widetilde{r}_{\Rm}(\cdot,t)$ and $\widetilde{r}_{\Ric}(\cdot,t)$ are $1$-Lipschitz continuous with respect to the metric $d_{g(t)}$. This is proven in the exact same way as the corresponding results for the parabolic scales $r_{\Rm}(\cdot,t)$ and $r_{\Ric}(\cdot,t)$ in Theorem \ref{thm.lipschitzholder} and Theorem \ref{thm.riccilipschitzholder}.

A straightforward consequence of the Lipschitz continuity is the following local Harnack type inequality: if $(p_0,t_0) \in M \times [0,T)$ is a point in a Ricci flow (say not identically Ricci flat), then
\begin{equation}\label{eq.bddscalarharnack}
\tfrac{1}{2}\, \widetilde{r}_{\Ric}(p_0,t_0)\le \widetilde{r}_{\Ric} \le \tfrac{3}{2}\, \widetilde{r}_{\Ric}(p_0,t_0), \quad \text{on} \, B_{g(t_0)}(p_0,\tfrac{1}{2}\, \widetilde{r}_{\Ric}(p_0,t_0)).
\end{equation}
To prove the integral concentration of the time-slice Ricci scale near singular points let then $(p_0,t_0)$ be a point satisfying the assumption of the lemma and set $r_0 := \widetilde{r}_{\Ric}(p_0,t_0)$. From the discussion following Definition \ref{def.noncollapsed}, we have that $\mu_{g(t_0)}(B_{g(t_0)}(p_0,\delta r_0)) \ge \kappa_1 \delta^{n}r_0^n$, where $\kappa_1=\kappa_1(n,i_0,k_0,T,R_0)$. Moreover, we can use (\ref{eq.bddscalarharnack}) to get $\widetilde{r}_{\Ric} \le \frac{3}{2}r_0$ on $B_{g(t_0)}(p_0,\delta r_0)$. Then we easily compute
\begin{equation*}
\int_{B_{g(t_0)}(p_0,\delta r_0)} \widetilde{r}_{\Ric}^{\,-2 \alpha} \, d\mu_{g(t_0)} \ge \kappa_1 \delta^n r_0^n (\tfrac{3}{2}\, r_0)^{-2 \alpha} \ge \kappa_1 \delta^n(\tfrac{2}{3})^{2 \alpha} \min \Big\{R_0^{\frac{2\alpha-n}{2}},1 \Big\} =:C_5>0.
\end{equation*}
In the second inequality we have used $n-2 \alpha \le 0$.
\end{proof}

We are now ready to prove the non-existence of well-behaved singularities in dimensions less than eight.
\begin{proof}[Proof of Theorem \ref{th.nosingularities}]
Without loss of generality $n \ge 4$. Fix any $t_1 \in (0,T)$ and $\delta \in (0,1)$.
Assume towards a contradiction that there exists a blow-up sequence $(p_i,t_i)$ along which $r_i:=r_{\Ric}(p_i,t_i) \to 0$ and which is $\delta$-well behaved. By the Pseudolocality Proposition $3.2$ in \cite{bam} we see that $\widetilde{r}_{\Rm}$ and $r_{\Rm}$ are comparable for bounded scalar curvature Ricci flows, that is, there exists a constant $C=C(n,R_0,i_0,k_0,T)$ such that $C\widetilde{r}_{\Rm} \leq r_{\Rm}$ for all points $(p,t)$ with $t$ large enough. Together with Theorem \ref{thm.squareroot}, this implies in particular that $\widetilde{r}^2_{\Ric} \geq r^2_{\Ric} \geq \delta a_3 r_{\Rm} \geq C\delta a_3 \widetilde{r}_{\Rm}$ on the ball $B_{g(t_i)}(p_i,\sqrt{\delta} \widetilde{r}_{\Ric}(p_i,t_i))$ for $i$ large enough so that $t_i \ge t_1$.

Since we have assumed that $M$ is closed, the initial metric satisfies $\inj(M,g(0))>i_0$ and $\Ric_{g(0)}\ge -(n-1)k_0 g(0)$, for some $i_0>0$, $k_0>0$. We can then apply Corollary \ref{cor.theorem17bamler} with $\e <1/8$ to obtain for every $i$ and every $r'_i \in (0,1)$
\begin{equation}
\int_{B_{g(t_i)}(p_i,r'_i)} \widetilde{r}_{\Rm}^{\,-\alpha} \, d\mu_{g(t_i)} \le F (r_i')^{n-4+2\e},
\end{equation}
where $\alpha:=4-2\e$ and $F=F(n,R_0,i_0,k_0,\e)$. In order to apply Lemma \ref{lemma.L8lowerbound}, notice that $\alpha \ge \frac{n}{2}$ if and only if $n<8$. Choosing $r_i':=\sqrt{\delta} r_i$ we can employ Lemma \ref{lemma.L8lowerbound}, and together with the inequality $\widetilde{r}_{\Ric}^{\, 2} \ge C\delta a_3 \widetilde{r}_{\Rm}$, we see that
\begin{align*}
0<C_5 &\le \int_{B_{g(t_i)}(p_i,r'_i)} \widetilde{r}_{\Ric}^{\, -2\alpha} \, d\mu_{g(t_i)}\\
&\le (C\delta a_3)^{-\alpha} \int_{B_{g(t_i)}(p_i,r'_i)} \widetilde{r}_{\Rm}^{\, -\alpha} \, d\mu_{g(t_i)} \\
&\le (C\delta a_3)^{-\alpha} F (\sqrt{\delta} r_i)^{n-4+2\e}.
\end{align*}
It is sufficient to let $i$ go to infinity to get the desired contradiction.

To prove the second statement in the theorem, assume that the flow cannot be extended past time $T$ and pick a point $p \in \Sigma$. Then by Theorem \ref{thm.bridge} as well as Proposition \ref{prop.singRicciscale}, we see that for every $t_i \nearrow T$, the sequence $(p,t_i)$ is a blow-up sequence, which is $\delta$-well behaved because $M=G_\delta$, in contradiction with what we proved above.
\end{proof}

We have seen in Theorem \ref{thm.squareroot} that a point $p$ almost maximising the Ricci curvature is well behaved. By Shi's estimates, also points in a sufficiently small neighbourhood of $p$ are almost maximising, namely on a scale comparable to its Riemann scale. It would be interesting to see which scale one can reach, considering that any improvement of this scale might be used (with the same argument as above) to exclude the singularity formation in low dimensions. Under the additional assumption of an injectivity radius bound as in \eqref{eq.injbound}, we are able to reach a sufficiently large scale, as the following proof shows.

\begin{proof}[Proof of Corollary \ref{cor.nosingularities}]
We can again assume without loss of generality that $n \ge 4$. Suppose towards a contradiction that the flow cannot be extended past time $T<\infty$. Then by Sesum's result in \cite{ses}, there exists a sequence $(p_i,t_i)$ with $t_i \nearrow T$ such that $\abs{\Ric}(p_i,t_i) = \sup_{M\times[0,t_i]} \abs{\Ric} \to \infty$. Set $r_i:=\widetilde{r}_{\Ric}(p_i,t_i)$. Notice that $r_i \to 0$. 

For $i$ large enough, due to the injectivity radius bound \eqref{eq.injbound}, we can apply Lemma~$1$ of Chen \cite{chen} (with $T=t_i$) to obtain the existence of some $\beta =\beta(\alpha)>0$ such that $\abs{\nabla\Ric} \leq \beta r_i^{-3}$. Setting $\delta=\frac{1}{2\beta}$, this yields 
\begin{equation*}
\abs{\Ric}(q,t_i) \geq \frac{1}{2} \abs{\Ric}(p_i,t_i) = \frac{1}{2} \sup_{M} \abs{\Ric}(\cdot,t_i)
\end{equation*}
for all $q\in B_{g(t_i)}(p_i,\delta r_i)$. For sufficiently large $i$ so that $r_i\le \min \{ R_0^{-1/2},\sqrt{t_i},1\}$, we can thus apply Part~$i)$ of Theorem \ref{thm.squareroot} to obtain $\widetilde{r}_{\Ric}^{\, 2}(q,t_i) \ge a_3 \widetilde{r}_{\Rm}(q,t_i)$ on $B_{g(t_i)}(p_i,\delta r_i)$. We can therefore argue as in the proof of Theorem \ref{th.nosingularities} to obtain the desired contradiction.
\end{proof}

We finish this article with a proofs of Theorem \ref{th.codimension8}. It heavily relies on our localised version of Sesum's result from Theorem \ref{thm.bridge}, as well as the fact that the Ricci curvature blows up at least at a Type~I rate at any singular point.

\begin{proof}[Proof of Theorem \ref{th.codimension8}]
First of all, by Theorem \ref{thm.bridge}, $\Sigma=\Sigma^{\Ric}$. By Corollay \ref{RiccitypeIoscicor}, we furthermore know that $p\in\Sigma=\Sigma^{\Ric}$ if and only if $r_{\Ric}(p,t)<\sqrt{T-t}$ for all $t$. In particular, for any $t\in[0,T)$ we know that $\Sigma \subseteq \{p \mid r_{\Ric}(p,t)<\sqrt{T-t}\}$. We therefore have the inclusion
\begin{equation}
\Sigma_\delta = \Sigma\cap G_\delta \subseteq \set{p \mid r_{\Ric}(p,t) < \sqrt{T-t}} \cap G_\delta
\end{equation}
for every $t\in[0,T)$. Pick $t_1$ sufficiently large, such that $\sqrt{T-t} \le \min \set{R_0^{-1/2},\sqrt{t},1}$ for all $t\in [t_1,T)$. We can then appeal to Part $ii)$ of Theorem \ref{thm.squareroot} to get the inclusions
\begin{equation*}
\{p \mid r_{\Ric}(p,t) < \sqrt{T-t} \} \cap G_\delta \subseteq \{p \mid r_{\Rm}(p,t) < (a_3 \delta)^{-1} (T-t)\} \cap G_\delta.
\end{equation*}
By the Pseudolocality Proposition $3.2$ in \cite{bam} we see that $\widetilde{r}_{\Rm}$ and $r_{\Rm}$ are comparable for bounded scalar curvature Ricci flows, that is, there exists a constant $C=C(n,R_0,i_0,k_0,T)$ such that $C\widetilde{r}_{\Rm} \leq r_{\Rm}$ for all points $(p,t)$ with $t\in[t_1,T)$. Hence, we see that
\begin{equation}
\{p \mid r_{\Ric}(p,t) < \sqrt{T-t} \} \cap G_\delta \subseteq \{p \mid \widetilde{r}_{\Rm}(p,t) < (C a_3 \delta)^{-1} (T-t)\} \cap G_\delta
\end{equation}
for every $t \in [t_1,T)$. Increasing $t_1$ possibly even further, we can then also ensure that $t_1\geq t_0$ (where $t_0$ is the constant from Theorem \ref{th.proposition64bamler}) as well as $s:=2(C a_3 \delta)^{-1} (T-t)<1$. This allows us to apply Theorem \ref{th.proposition64bamler} with this $s<1$ and $r:=\frac{1}{2}<1$, yielding for any exponent $d' \in (0,4)$ a constant $E'$ such that we obtain the upper bound for every $t \in [t_1,T)$
\begin{equation*}
\mu_{g(t)}(\Sigma_\delta \cap B_{g(t)}(p,\tfrac{1}{2})) \leq \mu_{g(t)}(\set{p \mid \widetilde{r}_{\Rm}(p,t)<sr} \cap B_{g(t)}(p,r))< E' s^{d'} r^n < E (\sqrt{T-t})^{2d'}.
\end{equation*}
Here we have set $E := E' 2^{d'}( Ca_3 \delta)^{-d'}$. We note that $d= 2d'\in(0,8)$. This concludes the proof.
\end{proof}

\makeatletter
\def\@listi{%
  \itemsep=0pt
  \parsep=1pt
  \topsep=1pt}
\makeatother
{\fontsize{10}{11}\selectfont

\printaddress

\end{document}